\documentclass[12pt,a4paper]{article}
\usepackage{latexsym,amssymb}
\usepackage{t1enc}

\usepackage{amsmath}
\usepackage{amsthm}
\usepackage[cp1250]{inputenc}
\usepackage[T1]{fontenc}
\usepackage[english]{babel}
\usepackage[dvips]{graphicx}
\usepackage[a4paper, left=3.0cm, right=3.0cm, top=2.5cm, bottom=3.4cm, headsep=1.2cm]{geometry}

\newtheorem{theorem}{Theorem}[section]
\newtheorem{lemma}[theorem]{Lemma}
\newtheorem{proposition}[theorem]{Proposition}
\newtheorem{corollary}[theorem]{Corollary}
\newtheorem{remark}[theorem]{Remark}

\def\R{\mathbb{R}}
\def\C{\mathbb{C}}
\def\N{\mathbb{N}}

\def\d{\mathrm{\,d}}
\def\lma{\lambda}
\def\part{\partial}

\def\Ind{\mathrm{Ind}}

\def\conv{\mathrm{conv}}

\begin{document}
\title{\textbf{{Forced periodic solutions for nonresonant parabolic equations on $\R^N$}}}
\author{ Aleksander \'{C}wiszewski (\footnote{The research supported by the National Science Center grant no. 2013/09/B/ST1/01963.} ),  Renata {\L}ukasiak (\footnote{Nicolaus Copernicus University, Faculty of Mathematics and Computer Science, ul. Chopina 12/18, 87-100 Toru\'{n}, Poland})}

\maketitle

\begin{abstract}
Criteria for the existence of $T$-periodic solutions of
nonautonomous parabolic equation $u_t = \Delta u + f(t,x,u)$,
$x\in\R^N$, $t>0$ with asymptotically linear $f$ will be provided.
It is expressed in terms of time average function $\widehat f$ of
the nonlinear term $f$ and the spectrum of the Laplace operator
$\Delta$ on $\R^N$. One of them says that if the derivative
$\widehat f_\infty$ of  $\widehat f$ at infinity does not interact
with the spectrum of $\Delta$, i.e. $\mathrm{Ker} (\Delta + \widehat
f_\infty)=\{0\}$, then the parabolic equation admits a $T$-periodic
solution. Another theorem is derived in the situation, where the
linearization at $0$ and infinity differ topologically, i.e. the
total multiplicities of positive eigenvalues of the averaged
linearizations at $0$ and $\infty$ are different mod $2$.
\end{abstract}


\section{Introduction}
We shall be concerned with time $T$-periodic solutions of the
following parabolic problem
\begin{equation}\label{23082013-0153}
\left\{   \begin{array}{l}
\displaystyle{\frac{\part u}{\part t}} (x,t) = \Delta u (x,t) + f(t,x,u(x,t)),  \, x\in\R^N,\, t>  0,\\
u(\cdot, t) \in H^1 (\R^N), \, t\geq 0,
\end{array}  \right.
\end{equation}
where  $\Delta$ is the Laplace operator (with respect to $x$) and a function $f:[0,+\infty)\times\R^N\times\R\to\R$ is
$T$-periodic in time:
\begin{equation}\label{23082013-1207}
f(t,x,u) = f(t+T,x,u) \ \mbox{ for all } \ t\geq 0, u\in\R \mbox{ and a.e. } x\in\R^N.
\end{equation}
Periodic problems for parabolic equations were widely studied by
many authors by use of various methods. Some early results are due
to Brezis and Nirenberg \cite{BreNir},  Amman and Zehnder
\cite{Amman-Zehnder}, Nkashama and Willem \cite{NkaWill}, Hirano
\cite{Hirano1, Hirano2}, Pr\"{u}ss \cite{Pruss}, Hess \cite{Hess},
Shioji \cite{Shioji} and many others; see also \cite{Vejvoda} and
the references therein. These results treat the case where
$\Omega$ is bounded and are based either on topological degree and
coincidence index techniques in the spaces of functions depending
both on $x$ and time $t$ or on the translation along trajectories
operator to which fixed point theory is applied. In this paper we
shall study the case $\Omega=\R^N$ by applying translation along
trajectories approach together with fixed point index and Henry's
averaging (see \cite{Henry}) as in \cite{Cwiszewski-CEJM} (for a
general reference see also \cite{Cwiszewski-hab}).  In this case the
semigroup compactness arguments are no longer valid (since the
Rellich-Kondrachov theorem on $\R^N$ does not hold and the semigroup
of bounded linear operators generated by the linear heat equation
$u_t = \Delta u$ on $\R^N$ is not compact). Therefore  adequate
topological fixed point theory for noncompact maps and the
adaptation of proper averaging techniques is required.

We shall assume that $f:[0,+\infty)\times\R^N\!\times\, \R \to\R$
is such that, for all \mbox{$t,s\in [0,+\infty)$}, $u,v\in \R$ and a.e. $x\in\R^N$,  one has
\begin{equation}\label{23082013-1208}
 f(t,\cdot, u) \mbox{ is measurable and } \   |f(t,x,0)| \leq m_0 (x);
\end{equation}
\begin{equation}\label{23082013-1100}
|f(t,x,u)-f(s,x,v)| \leq \left(\tilde k (x) + k(x)|u|\right) |t-s|^{\theta}+ l(s,x)|u-v|;
\end{equation}
\begin{equation}\label{23082013-1209}
\left(f(t,x,u)-f(t,x,v)\right) (u-v) \leq - a |u-v|^2 + b(x) |u-v|^2
\end{equation}
where $m_0\in L^2(\R^N)$, $\theta\in (0,1)$, $\tilde k \in
L^2(\R^N)$, $k=k_{0}+k_{\infty}$, $l=l_0+l_\infty$ with $k_{0},
l_0(t,\cdot)\!\in\! L^p (\R^N)$, $k_{\infty}, l_\infty(t,\cdot)
\!\in\! L^\infty (\R^N)$,  $\sup_{t\geq 0}
(\|l_0(t,\cdot)\|_{L^p}\! +\! \|l_\infty
(t,\cdot)\|_{L^\infty})\!<\! +\infty$, $a>0$ and $b\in L^p(\R ^N)$
where, if not stated otherwise,
$$
2<p<\infty \ \mbox{ for } \ N=1,2 \ \ \mbox{ and } \ \  N \leq p <\infty \ \mbox{ for } \ N\geq 3
$$
(see Remark \ref{03092014-1626} for examples of functions satisfying these assumptions).\\
\indent We shall consider the following linearization property of $f$ at zero
\begin{equation}\label{06092013-1208}
\lim_{u\to 0} \frac{f(t,x,u)}{u} = \alpha(t,x):= \alpha_0(t,x) -\alpha_\infty (t,x)
\end{equation}
and at infinity
\begin{equation}\label{03092013-0751}
\lim_{|u|\to \infty} \frac{f(t,x,u)}{u} = \omega (t,x) := \omega_0
(t,x) - \omega_\infty (t,x)
\end{equation}
for all $x\in\R^N$ and $t\geq 0$, where $\alpha_0 (t,\cdot), \omega_0 (t,\cdot) \in L^p (\R^N)$, $\alpha_\infty, \omega_\infty (t,\cdot)\in L^\infty (\R^N)$, $\alpha_\infty(t,x) \geq \bar \alpha_\infty$, $\omega_\infty (t,x) \geq \bar \omega_\infty$,  for all $t\geq 0$ and a.a. $x\in\R^N$, with some $\bar \alpha_\infty, \bar \omega_\infty >0$ and
$$
\sup_{t\geq 0} \left( \|\alpha_0 (t,\cdot)\|_{L^p}
+ \|\alpha_\infty (t,\cdot)\|_{L^\infty} + \| \omega_0 (t,\cdot)\|_{L^p}
+ \|\omega_\infty (t,\cdot)\|_{L^\infty}\right)<+\infty.
$$
We shall also assume that, for all $t,s \geq 0$ and a.e. $x \in \R^N$,
$$
\begin{array}{l}
|\alpha_0(t,x)-\alpha_0(s,x)|\leq k^0(x) |t-s|^\nu  \ \mbox{ and } \ |\alpha_\infty (t,x)-\alpha_\infty(s,x)| \leq k^\infty(x)|t-s|^\nu,\\
|\omega_0(t,x)-\omega_0(s,x)|\leq k^0(x)|t-s|^\nu  \ \mbox{ and } \ |\omega_\infty(t,x)-\omega_\infty(s,x)|\leq k^\infty(x)|t-s|^\nu
\end{array}
$$
with $k^0 \in L^p (\R^N)$, $k^\infty \in L^\infty (\R^N)$ and $\nu \in (0,1)$.

Our main results are the following criteria for the existence of $T$-periodic solutions.
\begin{theorem}\label{02092013-0753}
Suppose that $f$ satisfies conditions {\em (\ref{23082013-1207})},
{\em(\ref{23082013-1208})}, {\em(\ref{23082013-1100})},
{\em(\ref{23082013-1209})} and {\em (\ref{03092013-0751})}. If
\begin{equation} \label{11102013-0626-0}
\left\{
\begin{array}{l}
\displaystyle{\frac{\part u}{\part t}}(x,t) = \lma\Delta u(x,t) + \lma \omega (t,x) u(x,t), \ x\in\R^N, \, t>0,  \\
u(\cdot,t)\in H^1(\R^N),\  t\geq 0,
\end{array}
\right.
\end{equation}
has no nonzero $T$-periodic solutions, for $\lma\in (0,1]$ and
$\mathrm{Ker}\, (\Delta+\widehat\omega) = \{0 \}$, where $\widehat
\omega:\R^N\to\R$ is the time average function of $\omega$, given by
$\widehat \omega (x) := \frac{1}{T}\int_{0}^{T} \omega(t,x) \d t$,
then the equation {\em (\ref{23082013-0153})} admits a $T$-periodic
solution
$$u\in C([0,+\infty),
H^2(\R^N)) \cap C^1([0,+\infty),L^2(\R^N)).$$
\end{theorem}
\noindent Our second result applies in the case where there exists a
trivial periodic solution $u\equiv 0$ and the previous theorem does
not imply the existence of a nontrivial periodic solution.
\begin{theorem}\label{06092013-1209}
Suppose that all the assumptions of Theorem {\em
\ref{02092013-0753}} are satisfied and additionally that (\ref{06092013-1208}) holds. If the equation
\begin{equation}\label{03012014-1931-0}
\left\{
\begin{array}{l}
\displaystyle{\frac{\part u}{\part t}}(x,t) = \lma\Delta u(x,t) + \lma \alpha (t,x) u(x,t), \ x\in\R^N, \, t>0,  \\
u(\cdot,t)\in H^1(\R^N),\ t\geq 0,
\end{array}
\right.
\end{equation}
has no nonzero $T$-periodic solutions for $\lma \in (0,1]$ and
$\mathrm{Ker}\, (\Delta+\widehat\alpha) = \mathrm{Ker}\,
(\Delta+\widehat\omega) = \{0 \}$ and $m(0)\not\equiv m(\infty)$
$mod$ $2$, where $m (0)$ and $m (\infty)$ are the total
multiplicities of the positive eigenvalues of
$\Delta+\widehat{\alpha}$ and $\Delta+\widehat{\omega}$,
respectively, then the equation {\em (\ref{23082013-0153})} admits a
nontrivial $T$-periodic solution  $u\in  C([0,+\infty), H^2(\R^N))\cap
C^1([0,+\infty),L^2(\R^N))$.
\end{theorem}
\begin{remark} \label{03092014-1626} $\mbox{ }$\\{\em
\indent (a) Let us give an example of a class of functions
satisfying  (\ref{23082013-1208}), (\ref{23082013-1100}) and
(\ref{23082013-1209}). Consider $f:[0,+\infty) \times \R^N \times \R
\to \R$ given by
$$
f(t,x,u):= U(t,x) + V(t,x) u + g(W(t,x)u)
$$
with $U,V,W:[0,+\infty)\times \R^N\to \R$ such that $U(t,\cdot)\in
L^2(\R^N)$ for all $t\geq 0$, $\sup_{\tau \geq 0} \|U(\tau,
\cdot)\|_{L^2} < +\infty$, $V =V_0 +V_\infty$, $W=W_0+W_\infty$,
$V_0(t,\cdot), W_0(t,\cdot) \in L^p (\R^N)$  for all $t\geq 0$ and
$\sup_{\tau \geq 0} \left(\|V_0(\tau,\cdot)\|_{L^p} + \|
W_0(\tau,\cdot)\|_{L^p}\right)<\infty$ and $V_\infty, W_\infty\in
L^\infty ([0,+\infty)\times \R^N)$. Moreover we assume that  there
are $L_U\in L^2(\R^N)$, $L_V, L_W\in L^p(\R^N)+L^\infty(\R^N)$
such that, for all $t,s\geq 0$ and a.e. $x\in\R^N$,
$|U(t,x)-U(s,x)|\leq L_U (x)|t-s|^\theta$, $|V(t,x)-V(s,x)|\leq
L_V (x)|t-s|^\theta$ and $|W(t,x)-W(s,x)|\leq L_W(x)|t-s|^\theta$.
Furthermore, $g:\R\to \R$ is assumed to be a bounded Lipschitz
function with a constant $L>0$ such that $g(0)=0$ and $g'(0)$
exists. Then the assumptions (\ref{23082013-1208}) and
(\ref{23082013-1100}) are satisfied. If additionally there is
$a>0$ such that we have $V_\infty(t,x)+L|W_\infty(t,x)|\leq -a$
for all $t\geq 0$ and a.e. $x\in\R^N$, then (\ref{23082013-1209})
holds. As a concrete example one may give $f(t,x,u):= - 2a u +
\sin\left(au + b u (1+|x|)^{-\rho} |\cos t| \right)$ where $a,b>0$
and $\rho>1$ if $N=1,2$ and $\rho > N/p$ if $N\geq 3$. Moreover,
in this particular case $\omega_\infty \equiv 2a$, $\omega_0\equiv
0$, $\alpha_\infty \equiv a$,
$\alpha_0(t,x) = b (1+|x|)^{-\rho}|\cos t|$.\\
\indent (b) The appearance of terms from the space $L^p(\R^N)$ is essential
for our considerations. The function $\alpha_0$ (or  $\omega_0$)
assures that the positive part of the spectrum
$\sigma(\Delta+\widehat\alpha)$
(or $\sigma(\Delta+\widehat\omega)$) consists of a finite number of eigenvalues with finite dimensional eigenspaces (see Remark \ref{03092014-1732} for a more detailed discussion). That makes the numbers $m(0)$ and $m(\infty)$ well-defined, i.e. the formulation of the above result is correct. Note that in Theorem \ref{06092013-1209} the appearance of the nontrivial term either $\alpha_0$ or $\omega_0$ belonging to $L^p(\R^N)$ is necessary to satisfy the desired condition $m (0)\not\equiv m(\infty)$ mod $2$.\\
\indent (c) In this paper we focus on the case when there is no
resonance both at $0$ and at $\infty$, i.e. when the nonexistence
assumptions for  (\ref{11102013-0626-0}) and
(\ref{03012014-1931-0}) hold, respectively. From the technical
point of view they enable us to use continuation along the
parameter $\lambda$ up to $\lma=1$. For small parameter $\lma>0$
the lack of $T$-periodic assumptions of (\ref{11102013-0626-0})
and (\ref{03012014-1931-0}) is implied by the conditions
$\mathrm{Ker} (\Delta +\widehat\omega)=\{0\}$ and $\mathrm{Ker}
(\Delta+\widehat\alpha)=\{0\}$, respectively.  The lack of
nontrivial $T$-periodic solutions for the problem $\frac{\part
u}{\part t}=\lambda (\Delta u - \alpha_\infty(t,x)u + \alpha_0
(t,x) u)$, $u(\cdot,t)\in H^1(\R^N)$, $x\in \R^N$, $t>0$, is
obvious if $\alpha$ is independent of time (which is possible also
when $f$ depends on time). Moreover the nonexistence condition
also holds in the general case if, for instance
\begin{equation}\nonumber
\sup_{t \in [0,T]}\|\alpha_0(\cdot,t)\|_{L^p} < \left\{
\begin{array}{lcll}
\frac{p^{1/2p}\bar \alpha_\infty^{1-1/2p}}{2^{1/2p}},&&
\textnormal{ if } N=1,\ p>2,
\\
\frac{p^{1/p}\bar \alpha_\infty^{(1-1/p)}}{4^{1/p}},&& \textnormal{
if } N=2,\ p>2,
\\
\frac{\bar\alpha_\infty^{1-N/2p}}{(N/2p)^{N/2p} C(N)^{N/p}},&&
\textnormal{ if } N\geq 3,  N\leq p<\infty,
\end{array}  \right.
\end{equation}
\noindent where $C(N)>0$ is the constant in the Sobolev inequality
$\|u\|_{L^{\frac{2N}{N-2}}}\leq C(N) \|\nabla u\|_{L^2}$, $u\in
H^1(\R^N)$ (for details see Remark \ref{05092014-1148}).\\
\indent (d) The resonant case was considered in \cite{Cw-Luk-2015}. \hfill $\square$} \end{remark}

Following the tail estimates techniques of Wang \cite{Wang}, who studied attractors, and Prizzi
\cite{Prizzi-FM}, who studied stationary states and connecting
orbits by use of Conley index, we develop a fixed point index
setting applicable to parabolic equations on $\R^N$. We shall show
that the translation along trajectories operator  ${\bold
\Phi}_T:H^1(\R^N) \to H^1(\R^N)$ for (\ref{23082013-0153}) is
ultimately compact, i.e. belongs to the class of maps for which the
fixed point index $\Ind({\bold \Phi}_T, U)$, with respect to open
subsets of $H^1(\R^N)$, can be considered (see e.g.
\cite{Akhmerov-et-al}). Clearly the nontriviality of that index will
imply the existence of the fixed point of ${\bold \Phi}_T$ in $U$,
which is the starting point of the corresponding periodic solution.
In order to determine the index $\Ind ({\bold \Phi}_T, U)$, we use
an averaging method, i.e. we embed the equation
(\ref{23082013-0153}) into the family of problems
\begin{equation}\label{03012014-1930}
\left\{
\begin{array}{l}
\displaystyle{\frac{\part u}{\part t}} (x,t) =  \Delta u (x,t) +  f \left(\frac{t}{\lma},x,u(x,t)\right),  \, \,  x\in\R^N,\, t>0,\, \lma >0,\\
u(\cdot, t) \in H^1 (\R^N), \, t\geq 0.
\end{array}
\right.
\end{equation}
According to Henry's averaging principle solutions of
(\ref{03012014-1930}) converge as $\lma\to 0^+$ to a solution of the averaged equation
\begin{equation}\label{23082013-0155}
\left\{   \begin{array}{l}
\displaystyle{\frac{\part u}{\part t}} (x,t) = \Delta u (x,t) + \widehat f(x,u(x,t)), \ x\in\R^N, \, t>0,\\
u(\cdot, t) \in H^1 (\R^N), \, t>0,
\end{array}   \right.
\end{equation}
where the time average function $\widehat f:\R^N\times \R\to \R$ of
$f$ is given by
$$
\widehat f(x,u):=\frac{1}{T}\int_{0}^{T} f(t,x,u)\d t.
$$
Exploiting the tail estimate technique of Wang and
Prizzi together with an extension of Henry's averaging principle we
prove that asymptotic assumptions on $f$ imply a sort of a priori
bounds conditions, i.e. that there are no $\lma T$-periodic solution
of (\ref{03012014-1930}), for $\lma\in (0,1]$, with initial states
of large $H^1$ norm (in case of Theorem \ref{02092013-0753}) and
also of small $H^1$ norm (in case of Theorem \ref{06092013-1209}),
i.e. initial states of $\lma T$-periodic solutions are located
outside some open bounded set $U\subset H^1(\R^N)$. This enables us
to use a sort of the averaging index formula stating that
\begin{equation}\label{03012014-2346}
\Ind ({\bold \Phi}_T, U) = \lim_{t\to 0^+} \Ind (\widehat {\bold
\Phi}_t, U)
\end{equation}
where $\widehat {\bold \Phi}_t$ is the translation along
trajectories operator for (\ref{23082013-0155}). In computation of
$\Ind (\widehat {\bold \Phi}_t, U)$, for small $t>0$, the spectral
properties of the operators $\Delta + \widehat\alpha$ and $\Delta
+ \widehat \omega$ are crucial.  We strongly use the fact that
their essential spectrum is contained in $(-\infty,0)$ and the rest
consists of positive eigenvalues with finite dimensional eigenspaces
and that numbers $m(0)$ and $m(\infty)$ are well-defined (i.e.
finite).

The paper is organized as follows. In Section 2 we recall the
concept of ultimately compact maps and fixed point index theory. In
Section 3 we strengthen in a general setting of sectorial operators
the initial condition continuity property and Henry's averaging
principle. Section 4 is devoted the ultimate compactness property of
the translation operator. In Section 5 we adapt the ideas of
\cite{Cwiszewski-CEJM} to the case $\Omega=\R^N$, proving the
averaging index formula (\ref{03012014-2346}) as well as verify a
priori bounds conditions for $\lma T$-periodic solutions of
(\ref{03012014-1930}) with $\lma \in (0,1]$. Finally, in Section 6
the main results are proved.

\section{Preliminaries}

\noindent {\bf Notation}. If $X$ is a normed space with the norm
$\|\cdot \|$, then, for $x_0\in X$ and $r>0$, we  put
$B_X(x_0,r):=\{x\in X \mid \|x-x_0\|<r\}$. By $\part U$ and
$\overline{U}$ we denote the boundary and the closure of $U\subset
X$. $\conv\, V$ and $\overline{\conv}^{X}\, V$ stand for  the convex
hull and the closed (in $X$) convex hull of $V\subset X$,
respectively.
By $(\cdot,\cdot)_0$ is denoted the inner product in $X$.\\

\noindent {\bf Measure of noncompactness}. If $X$ is a Banach space
and $V\subset X$ is bounded, then by $\beta_X (V)$ we denote the
infimum over all $r>0$ such that $V$ can be covered with a finite
number of open balls of radius $r$.  Clearly $\beta_X (V)$ is finite
and it is called the Hausdorff measure of noncompactness of the set
$V$ in the space $X$. It is not hard to show that $\beta_X (V)=0$
implies that
$V$ is relatively compact in $X$. More properties of the measure of noncompactness can be found in \cite{Deimling} or \cite{Akhmerov-et-al}.\\

\noindent {\bf Fixed point index}. Below we recall basic definitions and facts from the fixed point index theory for ultimately compact maps. For details we refer to \cite{Akhmerov-et-al}.\\
\indent We say that a map $\Phi: D \to X$, defined on a subset $D$
of a Banach space $X$ is {\em ultimately compact} if $V\subset X$ is
such that $\overline{\conv} \,\Phi(V\cap D)=V$, then $V$ is compact.
We shall say that an ultimately compact map $\Phi:\overline{U} \to
X$, defined on the closure of an open bounded set $U\subset X$, is
called {\em admissible} if $\Phi(u)\neq u$ for all $u\in \part U$.
By an {\em admissible homotopy} between two admissible maps $\Phi_0,
\Phi_1:\overline{U}\to X$ we mean a continuous map ${\bold
\Psi}:\overline{U} \times [0,1]\to X$ such that ${\bold \Psi}(\cdot,
0) = \Phi_0, \, \, {\bold \Psi}(\cdot, 1)=\Phi_1$, ${\bold \Psi}
(u,\mu)\neq u$ for all $u\in \part U$ and $\mu\in [0,1]$, and, for
any $V\subset X$, if $\overline{\conv} \,{\bold \Psi} ( (V\cap
\overline{U}) \times [0,1]) = V$, then $V$ is relatively compact.
$\Phi_0, \Phi_1$ are called homotopic then. A fixed point index for
ultimately compact maps was constructed in \cite[1.6.3 and
3.5.6]{Akhmerov-et-al}. Basic properties of the fixed point index
are collected in the following
\begin{proposition}\label{23082013-0220} $\mbox{ }$\\
{\em (i)   (existence)} If $\Ind(\Phi, U)\neq 0$, then there exists $u\in U$ such that $\Phi(u)=u$.\\
{\em (ii)  (additivity)} If $U_1, U_2 \subset U$ are open and $\Phi
(u)\neq u$ for all $u\in \overline{U\setminus (U_1\cup U_2)}$, then
$$\Ind(\Phi, U)=\Ind(\Phi, U_1) + \Ind(\Phi, U_2).$$
{\em (iii) (homotopy invariance)} If $\Phi_0, \Phi_1:\overline{U}\to
X$ are homotopic, then
$$\Ind(\Phi_0, U)=\Ind (\Phi_1, U).$$
{\em (iv)  (normalization)} Let $u_0\in X$ and $\Phi_{u_0}:\overline
U\to X$ be defined by $\Phi_{u_0}(u)=u_0$ for all $u\in
\overline{U}$. Then $\Ind (\Phi_{u_0},U)$ is equal $0$ if
$u_0\not\in U$ and $1$ if $u_0\in U$.
\end{proposition}
\begin{remark} {\em  If $\Phi: \overline{U} \to X$ is a compact map then
$\Ind(\Phi, U)$ is equal to the Leray-Schauder index $\Ind_{LS}
(\Phi, U)$ (see e.g. \cite{Granas}).}
\end{remark}

\section{Remarks on abstract continuity and averaging principle}
Let $A:D(A)\to X$ be a sectorial operator such that for some $a>0$,
$A+a I$ has its spectrum  in the half-plane $\{z\in \C \mid
\mathrm{Re}\, z >0 \}$. Let  $X^\alpha$, $0\leq \alpha<1$, be the
fractional power space determined by $A+aI$. It is well-known that
there exist $C_0, C_\alpha>0$ such that for all $t>0$
$$
\| e^{-t A}u\|_\alpha\leq C_0 e^{at} \|u\|_\alpha \  \mbox{ for all } \ u\in X^\alpha,
$$
$$
 \|e^{-tA} u\|_{\alpha} \leq C_\alpha t^{-\alpha}e^{at} \|u\|_{0} \ \mbox{ for all } \ u\in X
$$
where $\{e^{-tA}\}_{t\geq 0}$ is the semigroup generated by $-A$.
Consider the equation
\begin{equation}\label{28012014-1055}
\left\{   \begin{array}{l}
\dot u(t) = - A u(t) + F(t,u(t)), \ t>0,\\
u(0) = \bar u,
\end{array}  \right.
\end{equation}
where $\bar u \in X^\alpha$ and $F:[0,+\infty) \times X^\alpha \to
X$ is such that  there exists $C\geq 0$ with
\begin{equation}\label{2356-20150921}
\|F(t,u)\|\leq C(1+\|u\|_\alpha) \mbox{ for all } u\in X, \, t>0,
\end{equation}
and, for any bounded $V\subset X^\alpha\times [0,+\infty)$ there exist $D, L\geq 0$  and $\theta\in (0,1)$ with the property
\begin{equation}\label{2357-20150921}
\|F(t,u)-F(s,v)\|\leq D |t-s|^\theta + L \|u-v\|_\alpha \mbox{ for all } (u,t),(v,s) \in V.
\end{equation}
We shall say that $u: [0,+\infty) \to X^\alpha$ is a solution of above initial value
problem if
$$u\in C([0,+\infty),X^\alpha)\cap C((0,+\infty), D(A))\cap C^1((0,+\infty), X)
$$
and satisfies $(\ref{28012014-1055})$. By classical results (see
\cite{Cholewa} or \cite{Henry}), the problem $(\ref{28012014-1055})$
admits a unique global solution $u\in C([0,+\infty),X^\alpha)\cap
C((0,+\infty), D(A))\cap C^1((0,+\infty), X)$. Moreover, it is known
that $u$ being solution of $(\ref{28012014-1055})$ satisfies the
following Duhamel formula
\begin{equation}\label{22092015-1244}
u(t)= e^{-tA} u(0) + \int_0^t e^{-(t-s)A}F(s,u(s)) \d s,\qquad t>0.
\end{equation}
\begin{remark}\label{22092015-1343}{\em
Assume that $u:[0,T]\to X^\alpha$ is a solution of (\ref{28012014-1055}) with $T>0$. Then clearly, by (\ref{2356-20150921}) and (\ref{22092015-1244}), there is a constant $\tilde C=\tilde C(C,C_0, C_\alpha, a,T)>0$ such that for all $t\in (0,T]$
\begin{eqnarray*}
\|u(t)\|_\alpha & \leq & C_0 e^{at}\|\bar u\|_\alpha+\int_{0}^{t} C_\alpha (t-s)^{-\alpha} e^{a(t-s)}\|F(s,u(s))\|\d s\\
& \leq & \tilde C(1+\|\bar u\|_\alpha) + \tilde C \int_{0}^{t}(t-s)^{-\alpha} \|u(s)\|_\alpha\d s.
\end{eqnarray*}
This in view of \cite[Lemma 1.2.9]{Cholewa} implies that there exists $\bar C=\bar C(C,C_0, C_\alpha, a,T,\alpha)>0$ such that
\begin{equation}
\|u(t)\|_\alpha \leq \bar C(1+\|\bar u\|_\alpha) \ \mbox{ for all }\  t\in [0,T].
\end{equation}\nonumber
}\end{remark}
\begin{theorem}\label{30082013-0308}
Assume that $(\alpha_n)$ is a sequence of positive numbers such that $\alpha_n\to \alpha_0$ as $n\to +\infty$ for some $\alpha_0>0$ and that $A_n:= \alpha_n A$ for $n\geq 0$.
Let $F_n : [0,T]\times X^\alpha\to X$, $T>0$ , $n\geq 0$, satisfy {\em (\ref{2356-20150921})} and {\em(\ref{2357-20150921})} with common constants $C, L$ {\em(}independent of $n${\em)} and let, for each $u\in
X^\alpha$,
$$
\int_{0}^{t} F_n (s,u) \d s\to \int_{0}^{t} F_0(s,u) \d s  \ \
\mbox{ in } \ \ X  \mbox{ as } n\to +\infty
$$
uniformly with respect to $t\in [0,T]$. If
$u_n:[0,T]\to X^\alpha$, $n\geq 0$, are solutions of
$$
\dot  u(t)=-A_n u(t)+F_n (t,u(t)), \ t\in [0,T],
$$
and $u_n(0) \to u_0(0)$ in $X$, then $u_n (t)\to u_0 (t)$ in $X^\alpha$
uniformly with respect to $t$ from compact subsets of $(0,T]$.
\end{theorem}
\begin{remark} {\em
Recall that Henry's result from \cite{Henry} states that, under the
above assumptions with $\alpha_n \equiv 1$, if $u_n (0) \to u_0 (0)$ in $X^\alpha$, as
$n\to+\infty$, then $u_n(t)\to u_0(t)$ in $X^\alpha$ uniformly on
compact subsets of $[0, T)$. Here, inspired by the proof of
Proposition 2.3 of \cite{Prizzi-FM}, we modify Henry's proof.}
\end{remark}
\noindent In the proof we shall use the following lemma.
\begin{lemma}\label{30082013-2347}
Under the assumptions of Theorem {\em \ref{30082013-0308}}, for any
continuous  $u:[0,T]\to X^\alpha$,
$$
\int_{0}^{t} e^{-(t-s)A_n} F_n (s, u(s)) \d s \to \int_{0}^{t} e^{-(t-s)A_0} F_0 (s, u(s)) \d s \ \mbox{ in } X^\alpha \mbox{ as } \ n
\to +\infty,
$$
uniformly with respect to $t\in [0,T]$.
\end{lemma}
\noindent  {\bf Proof:} We shall adjust arguments from the proof of \cite[Lemma 3.4.7]{Henry}. First observe that due to the assumptions concerning the constant $L$ for $F_n$'s it is sufficient to show the assertion for $u\equiv \bar u$ where $\bar u\in X^{\alpha}$.
Take any $\varepsilon>0$. There exist $\delta>0$, $\tilde C>0$ and $\tilde a>0$ such that, for any $n\geq 0$ and $t\in [0,\delta]$
\begin{eqnarray}\nonumber
\left\|  \int_{0}^{t} e^{ - (t-s) A_n } F_n (s,\bar u) \d s \right\|_\alpha \leq
\int_{0}^{t} C_\alpha  \alpha_{n}^{-\alpha} (t-s)^{-\alpha} e^{\alpha_n a (t-s)}C (1+\|\bar u\|_\alpha) \d s\\
\leq \tilde C \int_{0}^{t} \tau^{-\alpha} e^{\tilde a \tau } \d \tau \leq  \tilde C e^{\tilde a T} (1-\alpha)^{-1}\delta ^{1-\alpha}\leq
 \tilde C e^{\tilde a T} \delta ^{1-\alpha}< \varepsilon/4 \label{2257-19102015}
\end{eqnarray}
and, for any $n\geq 0$ and $t  \in [\delta ,T]$,
\begin{equation}\label{2256-19102015}
\left\|  \int_{t-\delta}^{t} e^{- (t-s) A_n} F_n (s,\bar u) \d s \right\|_\alpha  \leq \tilde C \int_{0}^{\delta} \tau^{-\alpha} e^{\tilde a \tau } \d \tau \leq \tilde C e^{\tilde a T} (1-\alpha)^{-1}\delta ^{1-\alpha} < \varepsilon/4.
\end{equation}
Observe that, for any $n\geq 0$ and $t\in [\delta,T]$,
\begin{eqnarray*}
 \int_{0}^{t-\delta} e^{-(t-s)A_n} F_n (s, \bar u) \d s  = e^{-tA_n}\int_{0}^{t} F_n (\tau,\bar u)\d\tau-e^{-\delta A_n}\int_{t-\delta}^{t} F_n (\tau, \bar u)\d \tau\\
 + \int_{0}^{t-\delta} A_n e^{-(t-s)A_n}  \int_{s}^{t} F_n (\tau,\bar u) \d \tau\, \d s.
\end{eqnarray*}
Clearly,
$$
e^{-tA_n} \int_{0}^{t} F_n (\tau,\bar u)\d\tau \to e^{-t A_0 }\int_{0}^{t} F_0(\tau,\bar u)\d \tau, \  \mbox{ in } X^{\alpha},
$$
uniformly with respect to $t\in [\delta,T]$. Note also that, for all $t\in [\delta, T]$ and all $n\geq 1$,
$$
\left\|e^{-\delta A_n}\int_{t-\delta}^{t} F_n (\tau, \bar u)\d \tau\right\|_\alpha\leq  \tilde C e^{\tilde a \delta} \delta^{1-\alpha}\leq \varepsilon/4.
$$
Finally, for large $n$ and all $t\in [\delta, T]$ and $s\in [0,t-\delta]$, one has
\begin{eqnarray*}
& & \left\|A_n e^{-(t-s)A_n} \!\! \int_{s}^{t} F_n (\tau,\bar u) \d \tau \!-\! A_0 e^{-(t-s)A_0}\!\!  \int_{s}^{t} F_0 (\tau,\bar u) \d \tau \right\|_{\alpha}\\
& & \leq |\alpha_n-\alpha_0| \left\|A e^{-(t-s)A_n} \!\! \int_{s}^{t} F_n (\tau,\bar u) \d \tau \right\|_{\alpha} \\
& & \ \ \ \ \  + \alpha_0 \left\|A e^{-(t-s)A_n} \!\! \int_{s}^{t} F_n (\tau,\bar u) \d \tau \!-\! Ae^{-(t-s)A_0}\!\!  \int_{s}^{t} F_0 (\tau,\bar u) \d \tau \right\|_{\alpha} \\
& & \leq \bar C |\alpha_n-\alpha_0| \left\|e^{-(t-s)A_n} \!\! \int_{s}^{t} F_n (\tau,\bar u) \d \tau \right\|_{1+\alpha} \\
& & \ \ \ \ \  + \bar C \alpha_0 \left\|e^{-(t-s)A_n}\left(\int_{s}^{t} F_n (\tau,\bar u) \d \tau-\int_{s}^{t} F_0 (\tau,\bar u) \d \tau \right) \right\|_{1+\alpha}\\
& & \ \ \ \ \  + \bar C\alpha_0 \left\|\left(e^{-(t-s)A_n}-e^{-(t-s)A_0}\right)  \int_{s}^{t} F_0 (\tau,\bar u) \d \tau \right\|_{1+\alpha} \\
& & \leq |\alpha_n-\alpha_0|\frac{\bar C C_{1+\alpha} e^{\tilde a T}}{\delta^{1+\alpha}}\left\|\int_{s}^{t} F_n (\tau,\bar u) \d \tau  \right\| + \alpha_0 \frac{\bar C C_{1+\alpha}e^{\tilde a T}}{\delta^{1+\alpha}} \left\|\int_{s}^{t} F_n (\tau,\bar u) \d \tau-\int_{s}^{t} F_0 (\tau,\bar u) \d \tau  \right\| \\
& & \ \ \ \ \  + \alpha_0 \frac{\bar C C_{1+\alpha} e^{\tilde a T}}{(\alpha_0\delta/2)^{1+\alpha}}\left\| \left(e^{-((t-s)\alpha_n/\alpha_0 - \delta/2)A_0}-e^{-(t-s-\delta/2 )A_0}\right)  \int_{s}^{t} F_0 (\tau,\bar u) \d \tau \right\|
\end{eqnarray*}
where $\bar C>0$ is such that $\|Aw\|_\alpha\leq \bar C \|w \|_{1+\alpha}$ for all $w\in X^{1+\alpha}$.
Therefore for large $n$ and all $t\in [\delta,T]$
$$
\left\| \int_{0}^{t-\delta} e^{-(t-s)A_n} F_n (s, \bar u) \d s - \int_{0}^{t-\delta} e^{-(t-s)A_0} F_0 (s, \bar u) \d s \right\|_\alpha \leq
\varepsilon/4 + \varepsilon/4 + \varepsilon/4 = 3\varepsilon/4,
$$
which together with (\ref{2257-19102015}) and (\ref{2256-19102015}) ends the proof. \hfill $\square$\\

\noindent {\bf Proof of Theorem \ref{30082013-0308}:}
By the Duhamel formula, for $t\in (0,T]$ and $n\geq 1$,
\begin{eqnarray*}
u_n(t)\!-\!u_0(t) & = & e^{-tA_n}u_n(0)-e^{-t A_0}u_0(0) +\\
& & + \int_{0}^{t} e^{-(t-s)A_n}F_n(s, u_0(s))-e^{-(t-s)A_0}F_0(s,u_0(s)))\d s\\
& & +\int_{0}^{t} e^{-(t-s)A_n}(F_n(s, u_n(s))-F_n(s,u_0(s)))\d s.
\end{eqnarray*}
This gives, for all $t\in (0,T]$ and $n\geq 1$,
$$
\|u_n(t)-u_0(t)\|_{\alpha}\leq \gamma_n(t) + C_\alpha L\int_{0}^{t}
e^{a\alpha_n (t-s)}(\alpha_n(t-s))^{-\alpha} \|u_n(s)-u_0(s)\|_\alpha \d s
$$
with
\begin{align}
\gamma_n(t)&:= \!\frac{C_\alpha e^{a \alpha_n t}}{(\alpha_n t)^{\alpha}} \|u_n(0)-u_0(0)\|_{0}+ \left\|(e^{-t A_n}-e^{-t A_0})u_0(0)\right\|_{\alpha}\nonumber \\ &+
\left\|\int_{0}^{t}\!\! \left(e^{-(t-s)A_n}\!F_n(s, u_0(s))\!-\!e^{-(t-s)A_0}\! F_0(s,u_0(s))\right)\!\!\d
s\right\|_{\alpha}.\nonumber
\end{align}
This means that there are $\tilde a>0$ and $\tilde C>0$ such that, for all $t\in (0,T]$ and $n\geq 1$,
$$
\|u_n(t)-u_0(t)\|_{\alpha}\leq \gamma_n(t) + \tilde C\int_{0}^{t}
e^{\tilde a(t-s)}(t-s)^{-\alpha} \|u_n(s)-u_0(s)\|_\alpha \d s.
$$
By use of Lemma 7.1.1 of \cite{Henry}, we get
$$
\|u_n(t)-u_0(t)\|_\alpha \leq \gamma_n(t) +K \int_{0}^{t}
(t-s)^{-\alpha}\gamma_n(s)\d s
$$
for some constant $K>0$. Now let us fix $t \in [0,T]$ and take an
arbitrary $\delta \in (0, t)$. Observe also that
\begin{eqnarray*}
\int_{0}^{t} (t-s)^{-\alpha}\gamma_n(s)\d s \leq
\frac{2^\alpha}{\delta^\alpha} \int_{0}^{t-\delta/2} \gamma_n(s)\d s
+ \int_{t-\delta/2}^{t} (t-s)^{-\alpha}\gamma_n(s)\d s\\
\leq \frac{2^\alpha}{\delta^\alpha} \int_{0}^{T} \gamma_n(s)\d s  +
\frac{(\delta/2)^{1-\alpha}}{1-\alpha}\cdot \sup_{s\in [\delta/2,T]}
\gamma_n(s).
\end{eqnarray*}
Since, in view of Lemma \ref{30082013-2347}, $\gamma_n(t)\to 0$
uniformly with respect to $t$ from compact subsets of $(0,T]$ and the
functions $\gamma_n$, $n\geq 1$, are estimated from above by an integrable
function %(of the form $t\mapsto C (t^{-\alpha} + t^{1-\alpha})$ with some constant $C>0$),
we infer, by the dominated convergence theorem, that $\|u_n(t)-u_0(t)\|_{\alpha}\to 0$ as $n \to +\infty$ uniformly with respect to $t\in [\delta,T]$. \hfill $\square$\\

The above theorem allows us to strengthen Henry's averaging
principle. We assume that mappings $F_n:[0,+\infty)\times
X^\alpha \to X$, $n\geq 1$, satisfy (\ref{2356-20150921}) and (\ref{2357-20150921}) with common constants $C,L$ (independent of $n$) and that there exists $\widehat F: X^\alpha\to X$ such that, for all $\bar u\in X^\alpha$,
\begin{equation}\label{01092013-1402}
\lim_{\tau\to+\infty,\, n\to +\infty} \frac{1}{\tau} \int_{0}^{\tau}
F_n (t,\bar u) \d t = \widehat F (\bar u) \ \ \mbox{ in } \ \  X.
\end{equation}
\begin{theorem}\label{30082013-0314}
Suppose $F_n$ and $\widehat F$ are as above, $\bar u_n \to \bar u_0$ in $X$, $\lma_n \to 0^+$ as $n\to+\infty$, and
$u_n: [0,+\infty) \to X^\alpha$, $n\geq 1$, are solutions of
$$\left\{
\begin{array}{l}
\dot u(t) = - A u(t) + F_n (t/\lma_n,u(t)), \ t>0,\\
u(0)=\bar u_n.
\end{array}\right.
$$
Then $u_n(t)\to \widehat u(t)$  in $X^\alpha$ uniformly with respect to $t$ from compact subsets of
$(0,+\infty)$ where $\widehat u:[0,+\infty)\to X^\alpha$ is the solution of
$$\left\{
\begin{array}{l}
\dot u(t)=-Au(t)+\widehat F(u(t)), \, t>0,\\
\ u(0)=\bar u.
\end{array}\right.
$$
\end{theorem}
\noindent {\bf Proof:} Let $\tilde F_n:=F_n(\cdot/\lma_n,\cdot)$ and
$\tilde F_0:=\widehat F$.  Observe that, using
(\ref{01092013-1402}), we get, for any $\bar u \in X^\alpha$ and
$t>0$,
$$
\int_{0}^{t} \tilde F_n (s,\bar u)\d s=  \lma_n \int_{0}^{t/\lma_n}
F_n(\rho,\bar u) \d \rho\to t\widehat F(\bar u) = \int_{0}^{t} \tilde F_0 (\bar u) \d s \ \
\mbox{ in } X, \mbox{ as }  n\to+\infty.
$$
Clearly, $\tilde F_n$, $n\geq 1$, and $\tilde F_0$ satisfy (\ref{2356-20150921}) and (\ref{2357-20150921}) with the common constants $C,L$. It can be easily verified that the convergence above is uniform with respect to $t$ from bounded subintervals of $[0,+\infty)$.   Now, an application of  Theorem \ref{30082013-0308} yields the assertion. \hfill $\square$
\begin{remark} $\mbox{ }$\\
{\em \indent   (a)  The above result is an improvement of the continuation theorem and the Henry averaging principle \cite[Th. 3.4.9]{Henry} to the case when initial values from $X^\alpha$ converge in the topology of $X$ (not $X^\alpha$). This will appear crucial when establishing the ultimate compactness property and verifying a priori estimates in the proofs of main results. We shall need to consider solutions in the phase space $X^{1/2}=H^1(\R^N)$ (cf. Remark \ref{03092014-1626}) while the compactness of sequences of initial values is possible with respect to the $L^2(\R^N)$ topology only.\\
\indent (b) An averaging principle for parabolic equations on $\R^N$ was also proved in \cite{Antoci-Prizzi} where time dependent coefficients of the elliptic operator were considered. Here we have provided a general abstract approach.\\
}\end{remark}

\section{Continuity, averaging and compactness for the parabolic equation}
We transform (\ref{23082013-0153}) into an abstract
evolution equation. To this end define an operator ${\bold A}:D({\bold A})\to
X$  in the space $X:=L^2(\R^N)$ by
$$
{\bold A} u:=-\sum_{i,j=1}^{N} a_{ij} \frac{\part^2 u}{\part x_j
\part x_i} , \mbox{ for } u\in D(\bold A):=H^2(\R^N),
$$
where $a_{ij}\in\R$, $i,j=1,\ldots, N$, are such that
$$
\sum_{i,j=1}^{N} a_{ij}\xi_i \xi_j > 0,
%\theta_0 |\xi|^2
\mbox{ for any } \xi\in \R^N,
$$
and $a_{ij}=a_{ji}$ for $i,j=1,\ldots,N$. It is well-known that
${\bold A}$ is a self-adjoint, positive and sectorial operator in
$L^2(\R^N)$.\\
\indent Suppose that $f$ is as in Section 1 and define ${\bold F}:[0,+\infty) \times H^1(\R^N)\to
L^2(\R^2)$ by $[{\bold F}(t,u)](x):=f(t,x,u(x))$ for a.e.
$x\in\R^N$.
\begin{lemma}\label{09012014-1600}
Under the above assumptions there are constants $D >0$,
depending only on $k$, $\tilde k$, $N$ and $p$, $L>0$, depending only   $l$, $N$ and $p$, and $C>0$, depending only
on $m_0, l, N$ and $p$, such
that, for all $t_1, t_2 \geq 0$ and $u_1, u_2\in H^1(\R^N)$,
$$
\|{\bold F}(t_1,u_1)-{\bold F}(t_2, u_2)\|_{L^2} \leq
D(1+\|u_1\|_{H^1})|t_1-t_2|^{\theta} + L \|u_1-u_2\|_{H^1}\
\mbox{ and }
$$
$$
\|{\bold F}(t,u)\|_{L^2} \leq C(1+\|u\|_{H^1}) \mbox{ for any }
t\geq 0  \mbox{ and } u\in H^1(\R^N).
$$
\end{lemma}
\noindent Before we pass to the proof of Lemma \ref{09012014-1600} we shall provide the following technical
lemma.
\begin{lemma}\label{28032015-2355}
There exist constants $C_1=C_1(N,p)>0$ and
$C_2=C_2(N,p)>0$ such that for any $u \in H^1(\R^N)$
\begin{equation}
\|u\|_{L^{2p/(p-1)}}\leq C_1\|u\|_{H^1}.\label{20052015-1829}
\end{equation}
and
\begin{equation}\label{21032015-1447}
\|u\|_{L^{2p/(p-2)}}\leq C_2 \|u\|_{H^1}.
\end{equation}
\end{lemma}
\begin{proof}[\bf{Proof.}]
 Take any $u \in H^1(\R^N)$. If $N=1$, then by use of the H\"{o}lder and interpolation inequalities together with the continuity of the embedding  $H^1(\R)$ into $L^\infty(\R)$, one gets
\begin{equation}
\|u\|_{L^{2p/(p-1)}}\leq \|u\|_{L^2}^{1-1/p}\|u\|_{L^\infty}^{1/p}
\leq C \|u\|_{H^1}\nonumber
\end{equation}
and
\begin{equation}
\|u\|_{L^{2p/(p-2)}}\leq \|u\|_{L^2}^{1-2/p}\|u\|_{L^\infty}^{2/p}
\leq  C\|u\|_{H^1},\nonumber \label{22032015-2318}
\end{equation}
\noindent where $C>0$ is such that $\|v\|_{L^\infty}\leq
C\|v\|_{H^1}$ for all $v \in H^1(\R)$. If $N=2$, then
\begin{equation}
\|u\|_{L^{2p/(p-1)}}\leq \|u\|_{L^2}^{1-2/p}\|u\|_{L^4}^{2/p} \leq
C\|u\|_{H^1},\nonumber
\end{equation}
\noindent where $C>0$ is the constant from the inequality
$\|v\|_{L^4}\leq C\|v\|_{H^1}$ for all $v \in H^1(\R^2)$. Similarly,
\begin{equation}
\|u\|_{L^{2p/(p-2)}}\leq
\|u\|_{L^2}^{1-2q/(pq-2p)}\|u\|_{L^q}^{2q/(pq-2p)} \leq
C\|u\|_{H^1}, \nonumber
\end{equation}
where $q$ is an arbitrary fixed number from $
(\frac{2p}{p-2},+\infty)$ and  $C>0$ is the constant coming from the
fact that $H^1(\R^2)$ embeds continuously into $L^s(\R^2)$ for any
$s \in [2,+\infty)$. Finally, if $N\geq 3$, then, by the same
techniques, we get
\begin{equation}
\|u\|_{L^{2p/(p-1)}}\leq
\|u\|_{L^2}^{1-N/2p}\|u\|_{L^{2N/(N-2)}}^{N/2p} \leq C
\|u\|_{H^1}\nonumber
\end{equation}
and
\begin{equation}
\|u\|_{L^{2p/(p-2)}}\leq
\|u\|_{L^2}^{1-N/p}\|u\|_{L^{2N/(N-2)}}^{N/p} \leq C
\|u\|_{H^1},\label{22032015-2320} \nonumber
\end{equation}
\noindent where $C>0$ is the constant in the Sobolev inequality
$\|v\|_{L^{2N/(N-2)}} \leq C\|v\|_{H^1}$ for all $v\in H^1(\R^N)$.
\end{proof}
\noindent {\bf Proof of Lemma \ref{09012014-1600}}: By use of
(\ref{23082013-1100}), the H\"{o}lder inequality and Lemma \ref{28032015-2355}, one finds constants $D=D(k,\tilde k, N, p)>0$  and $L=L(l,N,p)>0$ such that, for any
$t_1,t_2\geq 0$ and $u_1,u_2\in H^1(\R^N)$,
\begin{eqnarray*}
\|{\bold F}(t_1,u_1)-{\bold F}(t_2, u_2)\|_{L^2} \!\leq\!(\|\tilde
k\|_{L^2} + C_2\|k_{0}\|_{L^p}\|u_1\|_{H^1}+
\|k_{\infty}\|_{L^\infty}\|u_1\|_{L^2}  )|t_1-t_2|^{\theta}
 \\ + C_2\|l_0 (t_2,\cdot)\|_{L^p}\|u_1-u_2\|_{H^1} + \|l_\infty(t_2,\cdot)\|_{L^\infty}\|u_1-u_2\|_{L^2}\\
 \leq D(1+\|u_1\|_{H^1})|t_1-t_2|^{\theta} + L\|u_1-u_2\|_{H^1}.
\end{eqnarray*}
\indent Furthermore, by (\ref{23082013-1100}), one also has
$|f(t,x,u)|\leq |f(t,x,0)| + l(t,x)|u|$ for $t\geq 0$,\, $x\in\R^N,
\, u\in\R.$ This gives the existence of $C=C(m_0, l, N, p)>0$ such that
$$
\|{\bold F}(t,u)\|_{L^2} \leq \|m_0\|_{L^2} + C_2 \|l_0(t,\cdot)\|_{L^p} \|u\|_{H^1} + \|l_\infty(t,\cdot)\|_{L^\infty}\|u\|_{L^2} \leq   C (1+\|u\|_{H^1})
$$
for any $t\geq 0$  and $u\in H^1(\R^N)$. \hfill  $\square$\\

Consider now the evolutionary problem
\begin{equation}\label{11112013-1916}
\dot u(t) = - {\bold A}u(t) +{\bold F} (t,u(t)), \, t\geq 0, \ \
u(0) = \bar u\in H^1(\R^N).
\end{equation}
Due to Lemma \ref{09012014-1600} and standard results in theory of
abstract evolution equations (see \cite{Henry} or \cite{Cholewa}),
the problem (\ref{11112013-1916}) admits a unique global solution
$u\in C([0,+\infty),H^1(\R^N))$ $\cap\ \  C((0,+\infty),
H^2(\R^N))\cap C^1((0,+\infty), L^2(\R^N))$. We shall say that
$u:[0,T_0) \to H^1(\R^N)$, $T_0>0$, is a {\em solution} ({\em
$H^1$-solution}) of
$$
\left\{\begin{array}{cl} \displaystyle{\frac{\part u}{\part t}}
(x,t) = {\cal A} u(x,t)+f(t,x,u(x,t)), & x\in\R^N, \, t\in (0,
T_0),\\ u(x,0)=\bar u (x), &  x\in\R^N,
\end{array}\right.
$$
for some $\bar u \in H^1(\R^N)$, where ${\cal A} = \sum_{i,j=1}^{N}
a_{ij} \frac{\part^2 }{\part x_j
\part x_i}$, if
$$u\in C([0,+\infty),H^1(\R^N))\cap C((0,+\infty), H^2(\R^N))\cap C^1((0,+\infty), L^2(\R^N))
$$
and (\ref{11112013-1916}) holds. In this sense we have global in
time existence and uniqueness of solutions for the parabolic partial
differential equation.

The continuity of solutions properties are collected below.
\begin{proposition}\label{28082013-0033}   {\em (compare \cite[Prop.
2.3]{Prizzi-FM})} Assume that functions $f_n:
[0,+\infty) \times \R^N\times \R\to \R$, $n\geq 0$, satisfy the
assumptions {\em (\ref{23082013-1208})}  with common $m_0$ and {\em
(\ref{23082013-1100})} with common $l$ and that $f_n(t,x,u) \to
f_0(t,x,u)$, for all $t\geq 0$, $u\in \R$, a.e. $x\in\R^N$, and
$f_n (t,\cdot,0)\to f_0(t,\cdot,0)$ in $L^2(\R^N)$ for all $t\geq
0$. Suppose that $(\alpha_n)$ is a sequence of positive numbers such that $\alpha_n \to \alpha_0$, as $n\to +\infty$, for some $\alpha_0>0$.
Let $u_n:[0,T] \to H^1(\R^N)$, $n\geq 0$, be a solution of
$$
\displaystyle{\frac{\part u}{\part t}}
(x,t) = \alpha_n {\cal A} u(x,t)+f_n (t,x,u(x,t)), \,  x\in\R^N, \, t\in (0, T],
$$
such that, for some $R >0$, $\|u_n (t)\|_{H^1} \leq R$, for all $t\in[0,T]$ and $n\geq 0$. Then
$f_n(t,\cdot, u(\cdot))\to f_0(t,\cdot, u(\cdot))$ in $L^2(\R^N)$ for any $u\in H^1(\R^N)$ and $t\geq 0$ and\\
{\em (i)} if $u_n(0)\to u_0(0)$ in $L^2(\R^N)$ as $n\to \infty$, then $u_n(t)\to u(t)$ in $H^1(\R^N)$ for $t$ from compact subsets of $(0,T]$.\\
{\em (ii)} if $u_n(0)\to u_0(0)$ in $H^1(\R^N)$  as $n\to\infty$, then
$u_n(t)\to u_0(t)$ in $H^1(\R ^N)$ uniformly for $t\in[0,T]$.
\end{proposition}
\noindent {\bf Proof}: Define ${\bold F}_n:[0,+\infty)\times
H^1(\R^N)\to L^2(\R^N)$, $n\geq 0$, by $[{\bold
F}_n(t,u)](x):=f_n(t,x,u(x))$. Note that, in view of
(\ref{23082013-1100}), for any $t\geq 0$ and $u\in H^1(\R^N)$ and
a.e. $x\in\R^N$
\begin{eqnarray*}
|f_n(t,x,u(x))\!-\!f_0(t,x,u(x))|^2 \!\leq\!
2|f_n(t,x,0)-f_0(t,x,0)|^2\! \!+\! 4 |l(x,t)u|^2.
\end{eqnarray*}
Since, for any $t\geq 0$, $f_n(t,\cdot,0) \to f_0(t,\cdot,0)$ in
$L^2(\R^N)$ as $n\to +\infty$, the right hand side can be estimated
by an integrated function, which due to the Lebesgue dominated
convergence theorem implies ${\bold F}_n(t,u)\to {\bold F}_0 (t,u)$
in $L^2(\R^N)$ as $ n \to +\infty$. Moreover, by use of Lemma
\ref{09012014-1600}, we may pass to the limit under the integral to
get $\int_{0}^{t} {\bold F}_n (s, u)\d s \to \int_{0}^{t} {\bold
F}_0 (s,u)\d s  \ \mbox{ in } \ L^2(\R^N)$ for any $u\in H^1(\R^N)$ and $t\geq 0$. This
in view of Theorem \ref{30082013-0308} implies the assertion (i).
The assertion (ii) comes from the standard continuity theorem from
\cite{Henry}. \hfill $\square$\\

\noindent Let us also state an averaging principle.
\begin{proposition}\label{28082013-0034} Assume that functions $f_n:
[0,+\infty) \times \R^N\times \R\to \R$, $n\geq 0$,  satisfy the assumptions of Proposition {\em \ref{28082013-0033}}
and additionally {\em (\ref{23082013-1207})}.  Suppose that $\bar u_n \to \bar u_0$ in $L^2(\R^N)$, $\lma_n\to 0^+$ as $n\to +\infty$ and that $u_n: [0,+\infty) \to H^1(\R^N)$, $n\geq 1$, are solutions of
$$
\left\{\begin{array}{cl} \displaystyle{\frac{\part u}{\part t}} = {\cal A} u+f_n(t/\lma_n,x,u), & x\in\R^N, \, t>0,\\
u(x,0)=\bar u_n (x), &  x\in\R^N.
\end{array}\right.
$$
Then $u_n(t)\to \widehat u(t)$ in $H^1(\R^N)$ uniformly on compact subsets of $(0,+\infty)$, where $\widehat u:[0,+\infty)\to H^1(\R^N)$
is the solution of
$$
\left\{\begin{array}{cl} \displaystyle{\frac{\part u}{\part t}} = {\cal A} u+\widehat f_0 (x,u), & x\in\R^N, \, t>0,\\
u(x,0)=\bar u_0 (x), &  x\in\R^N,
\end{array}\right.
$$
with $\widehat f_0:\R^N \times \R\to \R$ given by $\widehat f_0 (x,u):=\frac{1}{T}\int_{0}^{T} f_0(t,x,u)\d t$ for all $u\in\R$ and
a.e. $x\in\R^N$.
\end{proposition}
\noindent {\bf Proof}: Define ${\bold F_n}:[0,+\infty)\times H^1(\R^N) \to L^2(\R^N)$, $n\geq 0$, by $[{\bold F}_n(t,u)](x):=f_n(t,x,u(x))$
and $\widehat {\bold F}_0: H^1(\R^N)\to L^2(\R^N)$ by $\widehat {\bold F}_0 (u): = \frac{1}{T}\int_{0}^{T} {\bold F}_0(t,u)\d t$. Clearly, for all $u\in H^1(\R^N)$, $\widehat {\bold F}(u)(x) = \widehat f(x,u(x))$ for a.e. $x\in\R^N$.
Fix any $\bar u\in H^1(\R^N)$ and $(\tau_n)$ in $(0,+\infty)$ such that $\tau_n\to +\infty$.
Clearly, ${\bold F}_n$, $n\geq 1$, are $T$-periodic in time. Consequently, one has
$$
I_n:=\frac{1}{\tau_n} \int_{0}^{\tau_n} {\bold F}_n (t,\bar u) \d t = \frac{[\tau_n/T]}{\tau_n/T} \cdot \frac{1}{T} \int_{0}^{T} {\bold F}_n (t,\bar u)\d t+
\frac{1}{\tau_n}\int_{0}^{\tau_n-[\tau_n/T]T} {\bold F}_n (t,\bar u)\d t.
$$
Hence, to see that $I_n\to \widehat {\bold F}_0(\bar u)$ it is is sufficient to prove that
$$
I_n^{(T)}:=\frac{1}{T}\int_{0}^{T} {\bold F}_n (t,\bar u)\d t \to \widehat {\bold F}_0 (\bar u) \ \mbox{ in } L^2(\R^N), \ \mbox{ as  } n\to+\infty.
$$
To this end observe that, for a.e. $x\in\R^N$,
$$
I_{n}^{(T)} (x) = \frac{1}{T} \int_{0}^{T} \!\! f_n(t,x,\bar u(x)) \d t\to  \frac{1}{T} \int_{0}^{T}\!\!  f_0 (t,x,\bar u(x)) \d t = \widehat f_0(x,u(x))=
[\widehat {\bold F}_0 (\bar u)](x).
$$
Moreover, by use of the assumptions on $f_n$'s, one has
\begin{eqnarray*}
|I_{n}^{(T)} (x)| & = & \left|\frac{1}{T} \int_{0}^{T} f_n (t,x,\bar u(x)) \d t\right| \leq m_0 (x)+ g(x)
\end{eqnarray*}
where  $g(x):= \frac{1}{T} \int_{0}^{T} |l (t,x)||\bar u(x)| \d t$.
and, by  use of Jensen's inequality,
$$
\int_{\R^N} |g(x)|^2 \d x \leq \frac{1}{T} \int_{\R^N}\int_{0}^{T} |l(t,x)|^2 |\bar u(x)|^2 \d t \d x <+\infty
$$
(see the proof of Lemma \ref{09012014-1600}). Hence, by the dominated convergence theorem we infer that
$I_n^{(T)}\to \widehat {\bold F}_0(\bar u)$ in $L^2(\R^N)$. Since $(\tau_n)$ was arbitrary it follows that
$$
\lim_{\tau \to+\infty,\ n\to +\infty} \frac{1}{\tau}\int_{0}^{\tau} {\bold F}_n (t,\bar u)\d t \to \widehat {\bold F}_0 (\bar u).
$$
Finally, we get the assertion by use of Theorem  \ref{30082013-0314}. \hfill $\square$\\

Now we pass to compactness issues that we treat with use of tail estimates technique.
\begin{lemma} \label{26082013-2204}
Assume that $f:[0,+\infty)\times \R^N\times \R\to\R$ satisfies {\em
(\ref{23082013-1208})}, {\em (\ref{23082013-1100})} and {\em
(\ref{23082013-1209})}. Suppose that $u:[0,T]\to H^{1}(\R^N)$ is a
solution of {\em (\ref{23082013-0153})} such that
$\|u(t)\|_{H^1}\leq R$ for all $t\in [0,T]$. Then there exists a
sequence $(\alpha_n)$ with $\alpha_n\to 0$ as $n\to\infty$ such that
$$
\int_{\R^N \setminus B(0,n)} |u(t)|^2 \d x \leq R^2 e^{-2 a t} +
\alpha_n \ \ \mbox{ for all } t\in [0,T], \, n\geq 1,
$$
where $\alpha_n's$ depend only on $N,p, R, m_0, a, b$ and  $a_{ij}'s$.
\end{lemma}
\noindent {\bf Proof}: it goes along the lines of \cite[Prop.
2.2]{Prizzi-FM}. The only difference is that here we have the
modified dissipativity condition (\ref{23082013-1209}), i.e.,
(\ref{23082013-1209}) implies
$$
f(t,x,u)u\leq - a|u|^2 + b(x)|u|^2  + f(t,x,0)u
$$
for $t\geq 0$, $x\in\R^N$, $u\in \R$, and one needs to modify the proof in a rather obvious way. \hfill $\square$\\

Now suppose that $a_{ij}\in C([0,1],\R)$, $i,j=1,\ldots, N$, are
such that $\sum_{i,j=1}^{N} a_{ij}(\mu)\xi_i \xi_j >0$ for any $\xi\in \R^N$ and $\mu\in [0,1]$. Let ${\bold
A}^{(\mu)}: D({\bold A}^{(\mu)}) \to L^2(\R^N)$, $\mu\in [0,1]$, be
given by
$$
{\bold A}^{(\mu)} u:= - \sum_{i,j=1}^{N} a_{ij}(\mu) \frac{\part^2
u}{\part x_j\part x_i}, \ u\in D({\bold A}^{(\mu)}):=H^2(\R^N).
$$
Let $h:[0,+\infty)\!\times\! \R^N\! \times\! \R\! \times\! [0,1]\to \R$ be such that, for all $t,s \geq 0$, $u,v\in
\R$, $\mu,\nu \in [0,1]$ and a.e. $x\in\R^N$,
\begin{equation}\label{16102013-0007}
h(t,\cdot,u,\mu) \mbox{ is measurable and }  |h(t,x,0,\mu)| \leq m_0(x),
\end{equation}
\begin{equation}\label{18032014-2500}
|h(t,x,u,\mu)-h(s,x,v,\mu)| \leq (\tilde k (x) +
k(x)|u|)|t-s|^{\theta} + l(s,x)|u-v|,
\end{equation}
\begin{equation}\label{16102013-0008}
|h(t,x,u,\mu) - h(t,x,u,\nu)|\leq l(t,x)\, |u|\,
|\rho(\mu)-\rho(\nu)|,
\end{equation}
\begin{equation}\label{16102013-0006}
(h(t,x,u,\mu)-h(t,x,v,\mu)) (u-v) \leq - a |u-v|^2 + b(x) |u-v|^2
\end{equation}
where  $m_0\in L^2(\R^N)$, $\theta\in (0,1)$, $\tilde k \in L^2(\R^N)$,
$k=k_{0}+k_{\infty}$ with $k_{0}\in L^p (\R^N)$, $k_{\infty}\in L^\infty
(\R^N)$, $l=l_0+l_\infty$ with $l_0 (t,\cdot)\in L^p(\R^N)$,
$l_\infty (t,\cdot)\in L^\infty(\R^N)$ for all $t\geq 0$, and
$\sup_{t\geq 0} \, (\|l_0 (t,\cdot)\|_{L^p} + \|l_\infty
(t,\cdot)\|_{L^\infty})<+\infty$,
$\rho\in C([0,1],\R)$, $a>0$ and $b\in L^p(\R ^N)$.\\ \ \\
\indent Under these assumptions consider
\begin{equation}\label{29082013-1258}
\dot u (t) = - {\bold A}^{(\mu)} u(t) + {\bold H}(t,u(t),\mu), \
t>0,
\end{equation}
where ${\bold H}:[0,+\infty) \times H^1(\R^N) \times [0,1]\to
L^2(\R^N)$ is defined by
$$
[{\bold H}(t,u,\mu)](x):= h(t,x, u(x),\mu) \mbox{ for } t\geq 0,\,
u\in H^{1}(\R^N),\, \mu\in[0,1], \mbox{ a.e. } x\in\R^N.
$$
Clearly, due to Lemma \ref{09012014-1600}, we get the existence and
uniqueness of solutions on $[0,+\infty)$.
Denote by $u(\cdot; \bar u, \mu)$ the solution of (\ref{29082013-1258}) satisfying the initial value condition $u(0)=\bar u$.\\
\indent The following tail estimates will be crucial in studying the
compactness properties of the translation along trajectories
operator of (\ref{29082013-1258}).
\begin{lemma}\label{26082013-2147}
Take any $\bar u_1, \bar u_2\in H^{1}(\R^N)$ and $\mu_1, \mu_2\in
[0,1]$ and suppose that there are solutions $u(\cdot;\bar
u_i,\mu_i):[0,T]\to H^1(\R^N)$, $i=1,2$ of {\em
(\ref{29082013-1258})}, for some  fixed $T>0$. If $\|u(t;\bar
u_1,\mu_1)\|_{H^1}\leq R$ and $\|u(t;\bar u_2,\mu_2)\|_{H^1} \leq R$
for all  $t\in [0,T]$ and some fixed $R>0$, then there exists a
sequence $(\alpha_n)$ with $\alpha_n\to 0$ as $n\to\infty$ such that
$$
\int_{\R^N\setminus B(0,n)} |u (t;\bar u_1,\mu_1)-u (t;\bar
u_2,\mu_2)|^2 \d x \leq e^{-2 a t} \| \bar u_1 - \bar u_2
\|_{L^2}^{2}  + Q \eta(\mu_1,\mu_2) + \alpha_n,
$$
for all $t\in [0,T]$ and $n\geq 1$, where $\alpha_n\geq 0$ and $Q>0$
depend only on $N, p, R, l, a, b$ and $a_{ij}'s$,
$$
\eta(\mu_1,\mu_2):= \max \left\{ |\rho(\mu_1)-\rho(\mu_2)|,
\max_{i,j=1,\ldots, N} |a_{ij} (\mu_1) - a_{ij} (\mu_2)| \right\}.
$$
\end{lemma}
\noindent {\bf Proof:} Let $\phi:[0,+\infty)\to \R$ be a smooth
function such that $\phi(s)\in [0,1]$ for $s\in [0,+\infty)$,
$\phi_{|[0,\frac{1}{2}]}\equiv 0$ and $\phi_{|[1,+\infty)}\equiv 1$
and let $\phi_n:\R^N\to\R$ be defined by $\phi_n
(x):=\phi(|x|^2/n^2)$, $x\in\R^N$. Put $u_1:=u(\cdot; \bar
u_1,\mu_1)$, $u_2:=u(\cdot; \bar u_2,\mu_2)$ and $v:=u_1-u_2$.
Observe that
\begin{eqnarray*}
\frac{1}{2}\frac{\d}{\d t} (v(t), \phi_n v(t))_{0} & =&
\frac{1}{2}\left(
 (v(t),\phi_n \dot v(t))_{0} + (\dot v(t), \phi_n v(t))_{0}\right) =  (\phi_n v(t),\dot v(t)))_{0} \\
& = & I_1(t)+I_2(t) + I_3(t)
\end{eqnarray*}
where
\begin{eqnarray*}
& I_1 (t) & := ( \phi_n v(t), - {\bold A}^{(\mu_1)}u_1(t) + {\bold A}^{(\mu_1)}u_2(t) )_0,\\
& I_2 (t) & := ( \phi_n v(t), -{\bold A}^{(\mu_1)}u_2(t) + {\bold A}^{(\mu_2)}u_2(t) )_0,\\
& I_3 (t) & := ( \phi_n v(t), {\bold H}  (t, u_1(t), \mu_1) - {\bold
H} (t, u_2(t), \mu_2) )_0.
\end{eqnarray*}
\indent As for the first term we notice that
\begin{eqnarray*}
I_1(t)& =&( \phi_n v(t), - {\bold A}^{(\mu_1)} v(t) )_0 \\
&   = &  -\int_{\R^N}\sum_{i,j=1}^{N} a_{ij}(\mu_1) \frac{\part
}{\part x_j}(\phi_n (x) v(t))  \frac{\part}{\part x_i}(v(t)) \d x\\
%&
%\end{eqnarray*}
%\begin{eqnarray*
& = & - \int_{\R^N} \phi_n (x) \sum_{i,j=1}^{N} a_{ij}(\mu_1)
\frac{\part }{\part x_j}( v(t))
\frac{\part}{\part x_i}(v(t)) \d x \\
& & -\frac{2}{n^2} \int_{\R^N} \sum_{i,j=1}^{N} \phi'(|x|^2/n^2)v(t) x_j a_{ij}(\mu_1) \frac{\part}{\part x_i} (v(t)) \d x\\
& \leq & \frac{2 L_\phi }{n^2} \int_{\big\{\frac{\sqrt{2}}{2}n\leq |x|\leq n\big\}} \sum_{i,j=1}^{N} a_{ij}(\mu_1) |x| |v(t)||\nabla_x v(t)|\d x\\
& \leq & \frac{2 L_\phi M N^2}{n} \|v(t)\|_{L^2}\|v(t)\|_{H^1}
\end{eqnarray*}
where $L_\phi:= \sup_{s\in[0,+\infty)}  |\phi'(s)|< \infty$ (as
$\phi'$ is smooth and nonzero on a bounded interval) and
$M:=\max_{1\leq i,j \leq N, \, \mu\in [0,1]} |a_{ij}(\mu)|$.
Further, in a similar manner
\begin{eqnarray*}
I_2(t) & = & - \int_{\R^N}\sum_{i,j=1}^{N}  \frac{\part }{ \part x_j} (\phi_n v(t)) (a_{ij}(\mu_1)-a_{ij}(\mu_2)) \frac{\part }{\part x_i }(u_2(t)) \d x  \\
& = & - \int_{\R^N} \phi_n (x) \sum_{i,j=1}^{N} (a_{ij}(\mu_1)- a_{ij}(\mu_2)) \frac{\part }{\part x_j}( v (t) ) \frac{\part}{\part x_i}( u_2(t) ) \d x\\
& & -\frac{2}{n^2} \int_{\R^N} \sum_{i,j=1}^{N} \phi'(|x|^2/n^2) v(t) x_j (a_{ij}(\mu_1)- a_{ij}(\mu_2))  \frac{\part}{\part x_i} (u_2 (t)) \d x\\
& \leq & \eta(\mu_1,\mu_2) \|v(t)\|_{H^1} \|u_2(t)\|_{H^1} + \frac{4
L_\phi  \eta(\mu_1,\mu_2) N^2}{n} \|v(t)\|_{L^2} \|u_2(t)\|_{H^1}.
\end{eqnarray*}
\indent To estimate $I_3(t)$  we see that (\ref{16102013-0006})
implies
\begin{eqnarray}
I_3(t)& = & \! \int_{\R^N}\!\! \phi_n(x) \left(
{\bold H}(t,u_1(t),\mu_1)\!-\! {\bold H}(t, u_2(t),\mu_2)\right) v(t) \d x \nonumber\\
& \leq &  \int_{\R^N} \phi_n(x) \left( {\bold H}(t,u_1(t),\mu_1)\!-\! {\bold H}(t, u_2(t),\mu_1)\right) v(t) \d x \nonumber\\
& & + \int_{\R^N} \phi_n(x) \left( {\bold H}(t,u_2(t),\mu_1)\!-\! {\bold H}(t, u_2(t),\mu_2)\right) v(t) \d x \nonumber \\
& \leq & - a\int_{\R^N}\phi_n(x)|v(t)|^2 \d x +
 \int_{\R^N}\!\!\!\phi_n(x)b(x)|v(t)|^2 \d x \nonumber \\
 & & +\, \int_{\R^N} l(t,x)|\rho(\mu_1)-\rho(\mu_2)||u_2(t)| |v(t)| \d x. \nonumber%\label{22012015-1604}
\end{eqnarray}
By use of the H\"{o}lder inequality together with Lemma
\ref{28032015-2355} one can get
\begin{eqnarray}
I_3(t)\!\! & & \leq \! - a\int_{\R^N}\phi_n(x)|v(t)|^2 \d x + C_1\|v(t)\|_{H^1}^2 \bigg(\int_{\big\{|x|\geq \frac{\sqrt{2}}{2}n\big\}}b(x)^p\d x\bigg)^{1/p} \nonumber \\
&  & \!\!\!\!\! + \eta(\mu_1,\mu_2)(C_2\|l_0(t,\cdot)\|_{L^p}\|u_2(t)\|_{H^1}
+\|l_{\infty}(t,\cdot)\|_{L^\infty}\|u_2(t)\|_{L^2})\|v(t)\|_{L^2}.\nonumber
\\\label{22012015-1512}
\end{eqnarray}
Hence we get, for any $n\geq 1$,
$$
\frac{\d}{\d t} (v(t), \phi_n v(t))_{0} \leq - 2 a (v(t), \phi_n
v(t))_{0} + \tilde C\eta(\mu_1,\mu_2) + \alpha_n
$$
for some constant $\tilde C=\tilde C(l, p, N, R)>0$, where
$(\alpha_n)_{n \in \mathbb{N}}$ is a sequence such that $\alpha_n
\to 0$ as $n \to +\infty$. Multiplying by $e^{2a t}$ and integrating
over $[0,\tau]$ one obtains
$$
e^{2a\tau}(v(\tau), \phi_n v(\tau))_{0}- (v(0), \phi_n v(0))_{0}
\leq  (2a)^{-1}(e^{2a\tau}-1)\, (\tilde C \eta (\mu_1,\mu_2)
+\alpha_n),
$$
which gives
$$
(v(\tau), \phi_n v(\tau))_{0} \leq e^{-2a \tau}\|v(0)\|^2_{L^2} +
(2a)^{-1} \left( \tilde C \eta (\mu_1,\mu_2)+\alpha_n\right).
$$
And this finally implies the assertion as
$\|\phi_n v(\tau)\|_{L^2}^{2}\leq (v(\tau), \phi_n v(\tau))_{0}$. \hfill $\square$\\

Let ${\bold \Psi}_t: H^1(\R^N)\times [0,1]\to H^1(\R^N)$, $t>0$, be
the translation operator for (\ref{29082013-1258}), i.e. ${\bold
\Psi}_t (\bar u,\mu) = u (t;\bar u,\mu)$ for $\bar u\in H^1(\R^N)$
and $\mu\in [0,1]$.
\begin{proposition}\label{01092013-1435}
Suppose that {\em (\ref{16102013-0007}),  (\ref{18032014-2500}), (\ref{16102013-0008})}   and {\em (\ref{16102013-0006})} are satisfied.\\
{\em (i)} For any bounded $V\subset H^1(\R^N)$ and $t>0$,
$\beta_{L^2} ({\bold \Psi}_t(V\times [0,1])) \leq e^{-at} \beta_{L^2}(V)$;\\
{\em (ii)} If a bounded $V\subset H^1(\R^N)$ is relatively compact
as a subset of $L^2(\R^N)$, then
${\bold \Psi}_t (V\times [0,1])$ is relatively compact in $H^1(\R^N)$;\\
{\em (iii)} If $V\subset \overline{\conv}^{H^1} {\bold
\Psi}_t(V\times [0,1])$ for some bounded $V\subset H^{1}(\R^N)$ and
$t>0$, then $V$ is relatively compact in $H^1(\R^N)$.
\end{proposition}
\noindent {\bf Proof:} (i) Observe that, for each $n\geq 1$,
$$
{\bold \Psi}_t(V\times [0,1])\subset \{u(t;\bar u, \mu)\mid \bar
u\in V, \, \mu\in [0,1]\}
 \subset W_n +R_n
$$
where $W_n:=\{\chi_n  u(t;\bar u,\mu) \mid \bar u\in V,\,\mu\in
[0,1]\}$ and $R_n:=\{ (1-\chi_n) u(t;\bar u,\mu) \mid \bar
u\in V,\,\mu\in [0,1]\}$ where $\chi_n$ is the characteristic
function of the ball $B(0,n)$. Note that $W_n$ may be viewed as a
subset of $H^1(B(0,n))$. Therefore, due to the Rellich-Kondrachov
theorem, $W_n$ is relatively compact in $L^2(\R^N)$. Hence
\begin{equation}\label{29082013-2327}
\beta_{L^2}({\bold \Psi}_t (V\times [0,1])) \leq \beta_{L^2}(R_n),
\ \ \mbox{ for all } \ \ n\geq 1.
\end{equation}
Now we need to estimate the measure of noncompactness of $R_n$ in
$L^2(\R^N)$. To this end fix an arbitrary $\varepsilon>0$. Choose a
finite covering of $V$ consisting of balls $B_{L^2} (\bar u_k,
r_\varepsilon)$, $k=1,\ldots, m_\varepsilon$, with
$r_\varepsilon:=\beta_{L^2} (V)+\varepsilon$ and such that $\bar u_k
\in V$ for each $k=1,\dots, m_{\varepsilon}$ and cover $[0,1]$ with
intervals  $(\mu_l-\delta, \mu_l+\delta)$, $l=1,\ldots, n_\delta$
where $\delta>0$ is such that $\eta(\mu_1, \mu_2)<\varepsilon$
whenever $|\mu_1-\mu_2|<\delta$.
  Put $\bar u_{k,l}:=(1-\chi_n) u(t;\bar u_k,\mu_l)$,
$k=1,\ldots,m_\varepsilon$, $l=1,\ldots, n_\delta$.\\
\indent Now take any $\bar v\in R_n$. There are $\bar u\in V$ and
$\mu\in [0,1]$ such that $\bar v = (1-\chi_n) u(t;\bar u,\mu )$.
Clearly there exist $k_0\in \{ 1,\ldots , m_\varepsilon \}$ and
$l_0\in \{1,\ldots, n_\delta \}$  such that $\|\bar u - \bar
u_{k_0}\| <r_\varepsilon$ and  $|\mu-\mu_{l_0}|<\delta.$ In view of
Lemma \ref{26082013-2147}
\begin{eqnarray*}
\|\bar v - \bar u_{k_0, l_0} \|_{L^2}^{2} & = &
\int_{\R^N\setminus B(0,n)} |u(t;\bar u,\mu)- u(t;\bar u_{k_0},\mu_{l_0} )|^2 \d x \\
&  \leq  & e^{-2 a t} \| \bar u-\bar u_{k_0} \|_{L^2}^{2}  + Q\, \eta(\mu,\mu_{l_0})  + \alpha_n \\
& \leq & r_{\varepsilon,n}:= e^{-2at}r_\varepsilon^2 +Q\,
\varepsilon  + \alpha_n,
\end{eqnarray*}
which means that $R_n$ is covered by the balls $B_{L^2}(\bar
u_{k,l},\sqrt{r_{\varepsilon,n}})$, $k=1,\ldots,m_\varepsilon$,
$l=1,\ldots, n_\delta$. This means that $\beta_{L^2}(R_n) \leq
\sqrt{r_{\varepsilon,n}}$ for any $\varepsilon>0$, and, in
consequence, $\beta_{L^2}(R_n)\leq (e^{-2at}(\beta_{L^2}(V))^2
+\alpha_n)^{1/2}$. Using (\ref{29082013-2327}) we get
$$
\beta_{L^2} ( {\bold \Psi}_t(V\times [0,1]) ) \leq
(e^{-2at}(\beta_{L^2}(V))^2 +\alpha_n)^{1/2},  \mbox{ for } n\geq 1.
$$
Finally, by a passage to the limit with $n\to \infty$ we obtain the required inequality as $\alpha_n \to 0^+$.\\
\indent (ii) Take any $(\bar u_n)$ in $V$ and $(\mu_n)$ in $[0,1]$.
We may assume that $\mu_n \to \mu_0$ for some $\mu_0 \in [0,1]$, as
$n \to +\infty$ . Since $(\bar u_n)$ is bounded, by the
Banach-Alaoglu theorem, we may suppose that $(\bar u_n)$ converges
weakly in $H^{1}(\R^N)$ to some $\bar u\in H^{1}(\R^N)$. By the
relative compactness of $V$ in $L^2(\R^N)$ we may assume that $\bar
u_n\to \bar u$ in $L^2(\R^N)$. Therefore, by use of Proposition
\ref{28082013-0033}, one has
${\bold \Psi}_t (\bar u_n, \mu_n) \to {\bold \Psi}_t(\bar u, \mu_0) $ in $H^1(\R^N)$, which ends the proof.\\
\indent (iii) Observe that here, by use of (i), one gets
$$
\beta_{L^2}(V) \leq \beta_{L^2} ({\bold \Psi}_t(V\times [0,1])) \leq
e^{-at}\beta_{L^2}(V).
$$
This implies $\beta_{L^2}(V)=0$, i.e. that $V$ is relatively compact
in $L^2(\R^N)$. To see that $V$ is relatively compact in $H^1(\R^N)$
observe that, by (ii), ${\bold \Psi}_t(V\times [0,1])$ is relatively
compact in $H^{1}(\R^N)$. \hfill $\square$

\section{Averaging index formula}

Consider the following parameterized equation
\begin{equation}\label{23102013-1303}
\left\{
\begin{array}{l}
\displaystyle{\frac{\part u}{\part t}} (x,t) = \Delta u (x,t) + h(t/\lma,x,u(x,t),\mu), \ t>0, \, x\in\R^N,\\
u(\cdot, t) \in H^1 (\R^N), \, t>0,
\end{array}  \right.
\end{equation}
where $h$ is as in the previous section and $\lambda>0$. Combining
the compactness result with averaging principle we get the following
result.
\begin{lemma}\label{23102013-1304}
Suppose $h$ satisfies conditions {\em (\ref{16102013-0007})}, {\em
(\ref{18032014-2500})}, {\em (\ref{16102013-0008})} and {\em
(\ref{16102013-0006})} and is $T$-periodic in the time variable
{\em($T>0$)}. If $(\bar u_n)$ is a bounded sequence in $H^1(\R^N)$,
$(\mu_n)$ in $[0,1]$, $(\lma_n)$ in $(0,+\infty)$ with $\lma_n \to
0^+$ as $n\to +\infty$ and $u_n: [0,+\infty) \to H^1(\R^N)$ are
solutions of {\em (\ref{23102013-1303})} with $\lma=\lma_n$,
$\mu=\mu_n$ such that $u_n(0)=u_n(\lma_n T)= \bar u_n$, then there
are a subsequence $(\bar u_{n_k})$ of $(\bar u_n)$ converging in
$H^1 (\R^N)$ to some $\bar u_0 \in H^2(\R^N)$ and a subsequence
$(\mu_{n_k})$ of $(\mu_n)$ converging to some $\mu_0\in [0,1]$, as
$k\to +\infty$, such that $\bar u_0$ is a solution of
$$
\Delta u(x)+\widehat h(x,u(x),\mu_0)=0, \ x\in \R^N,
$$
where $\widehat h:\!\R^N\! \times\! \R\times [0,1]\!\to\!\R$,
$\widehat h(x,u,\mu)\! :=
\displaystyle{\frac{1}{T}}\displaystyle{\int_{0}^{T}}\!\!
h(t,x,u,\mu) \d t$, $(x,u,\mu)\!\in\! \R^N\!\times \R\times [0,1]$.
Moreover, $u_{n_k} (t)\to \bar u_0$ in $H^1(\R^N)$, as $k\to
+\infty$,  uniformly with respect to $t$ from compact subsets of
$(0,+\infty)$.
\end{lemma}
\noindent {\bf Proof:} Recall that $u_n$ are solutions of $\dot u =
-{\bold A}u+{\bold H(t/\lma_n,  u,\mu_n)}$ with $u_n(0) = u_n(\lma_n
T)=\bar u_n$, $n\geq 1$, where ${\bold A}$ and ${\bold H}$ are as in
the previous section (with $a_{ij}=0$ if $i \neq j$ and $a_{ij}=1$ if $i=j$). Clearly, by the
sublinear growth, there exists $R>0$ such that $\|u_n(t)\|_{H^1}\leq
R$ for all $t>0$ and $n\geq 1$. For  an arbitrary
$M>0$ and $n\geq 1$ take $k_n\in \N$ such that $k_n \lma_n T>M$. In view of Lemma
\ref{26082013-2204}, for all  $m\geq 1$ and $n\geq 1$,
$$
\|(1-\chi_m)\bar u_n \|_{L^2}^2 = \|(1-\chi_m) u_n (k_n\lma_n
T)\|_{L^2}^2 \leq R^2 e^{-2a k_n \lma_n T} +\alpha_m \leq
R^2e^{-2aM}+\alpha_m,
$$
where $\chi_m$ is the characteristic function of $B(0,m)$. Since
$M>0$ is arbitrary we see that $\|(1-\chi_m) \bar u_n\|_{L^2}\leq
\sqrt{\alpha_m}$. Since, due to the Rellich-Kondrachov for any $m
\geq 1$, the set $\left\{\chi_m \bar u_n\right\}_{n\geq 1}$ is
relatively compact in $L^2(\R^N)$, we infer that $\{ \bar u_n
\}_{n\geq 1}$ is relatively compact in $L^2(\R^N)$. And since it is
bounded in $H^1(\R^N)$ we get a subsequence $(\bar u_{n_k})$,
denoted in the sequel again by $(\bar u_{n})$, such that $\bar
u_{n_k} \to \bar u_0$ in $L^2(\R^N)$ for some $\bar u_0\in
H^1(\R^N)$. We may also assume that $\mu_{n_k} \to \mu_0$ for some
$\mu_0\in [0,1]$. Hence, in view of Theorem \ref{30082013-0314},
$u_n(t)\to \widehat u(t)$ uniformly for $t$ from compact subsets of
$(0,+\infty)$ where $\widehat u:[0,+\infty)\to H^1(\R^N)$ is a
solution to
$$
\dot u = - {\bold A}u + \widehat {\bold H} (u,\mu_0), \, t>0,
$$
with $\widehat {\bold H}(u,\mu):= \frac{1}{T}\int_{0}^{T} {\bold
H}(t,u,\mu)\d t$ for $u\in H^1(\R^N)$, $\mu\in [0,1]$. Here note
that, for each $u\in H^1(\R^N)$ and $\mu\in [0,1]$,
$$
[\widehat {\bold H}(u,\mu)](x) = \widehat h(x,u(x),\mu) \  \mbox{
for all a.a. } x\in\R^N.
$$
Finally, for any $t>0$, we put $k_n:=[t/\lma_n T]$, $n\geq 1$,  and
see that
$$
\bar u_n = u_n(0)=u_n(k_n\lma_n T) \to \widehat u (t)  \ \mbox{ in }
H^1(\R^N), \mbox{ as } n\to+\infty.
$$
Hence $\widehat u(t)= \widehat u(0)=\bar u_0$ and $\bar u_n\to \bar
u_0$ in $H^1(\R^N)$. \hfill $\square$
\begin{remark}\label{15102013-0250}   {\em
Clearly that it follows from the proof of Lemma \ref{23102013-1304}
that
 if $f_n:[0,+\infty) \times \R^N\times \R \to \R$
 are as in Proposition \ref{28082013-0033} and satisfy (\ref{23082013-1209}) with common $a$ and $b$, then for any bounded sequence $(\bar u_n)$ in $H^1(\R^N)$, $(\lma_n)$ in $(0,+\infty)$ with   $\lma_n \to 0^+$ as $n\to +\infty$ and $u_n: [0,+\infty) \to H^1(\R^N)$ being $\lma_nT$-periodic solutions of
$$
\left\{\begin{array}{ll}
\displaystyle{\frac{\part u}{\part t}} = \Delta u + f_n(t/\lma_n,x, u), \, x\in\R^N, \, t>0,\\
u(x,0)= u(x,\lma_n T)=\bar u_n(x), \, x\in\R^N,
\end{array}\right.
$$
there is a subsequence $(\bar u_{n_k})$ of $(\bar u_n)$ converging
in $H^1 (\R^N)$ to some  $\bar u_0 \in H^2(\R^N)$  being a solution
of
$$
\Delta u(x)+\widehat f_0 (x,u(x))=0 \mbox{ on } \R^N.
$$
Moreover, $u_{n_k} (t)\to \bar u_0$ in $H^1(\R^N)$, as $k\to
+\infty$,  uniformly with respect to $t$ from compact subsets of
$(0,+\infty)$. \hfill $\square$}
\end{remark}

Now consider the following problem
\begin{equation}\label{23102013-0156}
\left\{
\begin{array}{l}
\displaystyle{\frac{\part u}{\part t}} (x,t) = \Delta u (x,t) + f\left(\frac{t}{\lma},x,u(x,t)\right), \ t>0, \, x\in\R^N,\\
u(\cdot, t) \in H^1 (\R^N), \, t>0,
\end{array}  \right.
\end{equation}
where $f$ satisfies conditions (\ref{23082013-1207}),
(\ref{23082013-1208}), (\ref{23082013-1100}) and
(\ref{23082013-1209}). We intend to prove an averaging index formula
that allows to express the fixed point index of translation along
trajectories operator for (\ref{23102013-0156}) in terms of the
averaged equation
\begin{equation} \label{23102013-0159}
\left\{
\begin{array}{l}
\displaystyle{\frac{\part u}{\part t}} (x,t) = \Delta u (x,t) + \widehat f(x,u(x,t)), \ t>0, \, x\in\R^N,\\
u(\cdot, t) \in H^1 (\R^N), \, t\geq 0,
\end{array}
\right.
\end{equation}
where  $\widehat f:\R^N\times\R\to\R$ is defined by
$$
\widehat f(x,u):= \frac{1}{T}\int_{0}^{T} f(t,x,u)\d t, \, x\in
\R^N, \, u\in\R.
$$
\begin{theorem}\label{26082013-2118}
Let $U\subset H^1(\R^N)$ be an open bounded set and by ${\bold
\Phi}_{t}^{(\lma)}$ and $\widehat {\bold \Phi}_t$, $t>0$, denote the
translation along trajectories operators (by time $t$) for the
equations {\em (\ref{23102013-0156})} and {\em
(\ref{23102013-0159})}, respectively. If the problem
\begin{equation}\label{23102013-01591}
\left\{
\begin{array}{l}
- \Delta u (x) =\widehat  f(x,u(x)), \, x\in\R^N,\\
u \in H^1 (\R^N),
\end{array}  \right.
\end{equation}
has no solution in $\part U$, then there exists $\lma_0>0$ such
that, for all $\lma\in(0,\lma_0]$, ${\bold \Phi}_{\lma
T}^{(\lma)}(\bar u)\neq \bar u$, $\widehat {\bold \Phi}_{\lma
T}(\bar u) \neq \bar u$ for all $\bar u \in \part U$, and
$$
\Ind ({\bold \Phi}_{\lma T}^{(\lma)}, U) = \Ind (\widehat {\bold
\Phi}_{\lma T}, U).
$$
\end{theorem}
\noindent {\bf Proof:} Define ${\bold H}:[0,+\infty) \times
H^1(\R^N)\times  [0,1]\to L^2(\R^N)$  by
$$
[{\bold H} (t,u,\mu)](x):= (1-\mu)  f(t,x,u(x)) + \mu \widehat f
(x,u(x)), \, \mbox{ for a.e. } x\in\R^N,
$$
and all $t>0, \, u\in H^1(\R^N)$. For a parameter $\lma>0$ consider
\begin{equation}\label{26082013-2052}
\dot u(t) = - {\bold A} u(t) + {\bold H} (t/\lma,u(t),\mu),\ t\in
[0,T],
\end{equation}
and the parameterized translation operator ${\bold
\Psi}_{t}^{(\lma)}:H^{1}(\R^N)\times [0,1]\to H^{1}(\R^N)$ defined
by
$$
{\bold \Psi}_{t}^{(\lma)} (\bar u, \mu):= u(t)
$$
where $u:[0,T]\to H^{1}(\R^N)$ is the solution of
(\ref{26082013-2052}) with $u(0)=\bar u$. Observe that for $\mu=0$,
(\ref{26082013-2052}) becomes
$$
\dot u (t) = -{\bold A} u(t)+ {\bold F}(t/\lma,u(t)), \ t\in [0,T],
$$
and we have ${\bold \Phi}_t^{(\lma)} = {\bold
\Psi}_{t}^{(\lma)}(\cdot,0)$. In the same way for $\mu=1$ the
equation (\ref{26082013-2052}) becomes
$$
\dot u (t)= -{\bold A} u(t) +\widehat {\bold F}(u(t)), \ t\in[0,T]
$$
and one has $\widehat {\bold \Phi}_t = {\bold \Psi}_t^{(\lma)}(\cdot, 1)$ (it does not depend on $\lma$).\\
\indent We claim that there exists $\lma_0>0$ such that, for all
$\lma\in (0,\lma_0]$,
\begin{equation}\label{26082013-2119}
{\bold \Psi}_{\lma T}^{(\lma)}(\bar u, \mu) \neq \bar u \ \ \mbox{
for all } \ \ \bar u\in \part U, \, \mu\in[0,1].
\end{equation}
Suppose the claim does not hold. Then there exist $(\bar u_n)$ in
$\part U$, $(\mu_n)$ in $[0,1]$ and $(\lma_n)$ with $\lma_n\to 0^+$
as $n\to\infty$ such that
$$
{\bold \Psi}_{\lma_n T}^{(\lma_n)}(\bar u_n, \mu_n) =\bar u_n \ \
\mbox{ for all } \ \ n\geq 1.
$$
This means that for each $n\geq 1$ there is a $\lma_n T$-periodic
solution $u_n:[0,+\infty) \to H^1 (\R^N)$ of (\ref{26082013-2052})
with $\lma=\lma_n$, $\mu=\mu_n$ and $u_n(0)=\bar u_n$. By Lemma
\ref{23102013-1304} we may assume that $\bar u_n \to \bar u_0$ in
$H^1(\R^N)$.
 Therefore $\bar u_0\in \part U\cap D({\bold A})$ and $0=-{\bold A}\bar u_0+\widehat {\bold F}(\bar u_0)$, a contradiction with the assumption. This proves the existence
of $\lma_0>0$ such that, for all $\lma \in (0,\lma_0]$, (\ref{26082013-2119}) holds.\\
\indent Now, due to Proposition \ref{01092013-1435} (iii), for each
$\lma\in (0,\lma_0]$, ${\bold \Psi}_{\lma T}^{(\lma)}$ is an
admissible homotopy in the sense of
fixed point index theory for ultimately compact maps. Finally, by Proposition \ref{23082013-0220}(iii), we get the desired equality of the indices. \hfill $\square$\\

As a consequence we get the following {\em continuation principle}.
\begin{corollary}\label{14102013-1246}
Suppose that  an open bounded $U\subset H^1(\R^N)$ is such that {\em
(\ref{23102013-01591})} has no solution in $\part U$, and for any
$\lma\in (0,1)$ the problem
\begin{equation}\label{09102013-1413}
\left\{
\begin{array}{l}
\displaystyle{\frac{\part u}{\part t}}  = \lma\Delta u +\lma f(t,x,u), x\in\R^N, \, t>0,\\
u(\cdot, t)\in H^1(\R^N), \, t\geq 0\\
u(x,0) = u(x,T), \, x\in \R^N,
\end{array}
\right.
\end{equation}
has no solution $u:[0,+\infty)\to H^1(\R^N)$ with $u(\cdot, 0)\in \part U$.
Then
$$
\Ind ({\bold \Phi}_T, U) = \lim_{t\to 0^+} \Ind ( \widehat {\bold
\Phi}_{t}, U)
$$
where ${\bold \Phi}_T$ is the translation along trajectories
operator for {\em (\ref{23082013-0153})}.
\end{corollary}
\noindent {\bf Proof:} Let $\lma_0>0$ be as in Theorem
\ref{26082013-2118}.  Since there are no solutions to
(\ref{09102013-1413}), we infer that
$$
{\bold \Phi}_{\lma T}^{(\lma)} (\bar u) \neq \bar u \mbox{ for any }
\bar u\in \part U,\, \lma\in (0,1).
$$
Now by Proposition \ref{01092013-1435} (iii) and the homotopy
invariance of the index, for any $\lambda \in (0,1]$, we get $\Ind
({\bold \Phi}_T, U)= \Ind(\widetilde {\bold \Phi}_T^{(1)}, U) = \Ind
(\widetilde {\bold \Phi}_{T}^{(\lma)}, U) = \Ind ({\bold \Phi}_{\lma
T}^{(\lma)}, U)$, where $\widetilde {\bold \Phi}_{T}^{(\lma)}$ is
the translation along trajectories operator for the parabolic
equation
in (\ref{09102013-1413}) with the parameter $\lambda$ and the last equality comes from  a time rescaling argument saying that $\widetilde {\bold \Phi}_{T}^{(\lma)} = {\bold \Phi}_{\lma T}^{(\lma)}$. Now an application of Theorem \ref{26082013-2118} completes the proof. \hfill $\square$\\

The rest of the section is devoted to methods of verification the
{\em a priori } bounds conditions  occurring in the above corollary
and computation of fixed point index. We shall use a linearization
approach.
\begin{proposition}\label{04102013-2013}
Suppose that $f$ satisfies conditions {\em (\ref{23082013-1207})}, {\em(\ref{23082013-1208})}, {\em (\ref{23082013-1100})}, {\em(\ref{23082013-1209})} and $f(t,x,0)=0$ for all $x\in \R^N$ and $t \geq 0$.\\
{\em (i)} \parbox[t]{140mm}{    If {\em (\ref{03092013-0751})}
holds, $\mathrm{Ker}\,(\Delta+\widehat \omega)=\{0\}$ and the linear
equation
\begin{equation}
\label{11102013-0626} \left\{
\begin{array}{l}
\displaystyle{\frac{\part u}{\part t}}(x,t) = \lma\Delta u(x,t) + \lma \omega(t,x) u(x,t), \ x\in\R^N, \, t>0,  \\
u(\cdot,t)\in H^1(\R^N), \ t\geq 0,
\end{array}
\right.
\end{equation}
has no nonzero $T$-periodic solutions for $\lma\in (0,1]$, then
there exists $R>0$ such that, for any $\lma\in (0,1]$
the problem {\em (\ref{09102013-1413})} has no $T$-periodic solutions $u:[0,+\infty) \to H^1(\R^N)$ with $\|u(0)\|_{H^1} \geq R$.}\\[1em]
{\em (ii)} \parbox[t]{140mm}{If {\em(\ref{06092013-1208})} holds,
$\mathrm{Ker}\,(\Delta+\widehat \alpha)=\{0\}$ and the linear
equation
\begin{equation}
\label{11102013-0627} \left\{
\begin{array}{l}
\displaystyle{\frac{\part u}{\part t}}(x,t) = \lma\Delta u(x,t) + \lma \alpha (t,x) u(x,t), \ x\in\R^N, \, t>0,  \\
u(\cdot,t)\in H^1(\R^N), \ t\geq 0,
\end{array}
\right.
\end{equation}
has no nonzero $T$-periodic solutions, then there exists $r>0$ such
that, for any $\lma\in (0,1]$ the problem {\em
(\ref{09102013-1413})} has no $T$-periodic solutions
$u:[0,+\infty)\to H^1(\R^N)$ with $0<\|u(0)\|_{H^1} \leq r$.}\\
[1em]
\end{proposition}
\noindent {\bf Proof:} (i) Suppose to the contrary, i.e. that for
any $n\geq 1$ there exist $\lma_n\in (0,1)$ and a time $T$-periodic
solution $u_n:[0, +\infty)\to H^1(\R^N)$ of
$$
\frac{\part u}{\part t} = \lma_n \Delta u + \lma_n f(t,x,u), \,
x\in\R^N, \, t>0
$$
with $\| u_n (0)\|_{H^1}\to +\infty$. This means that  $z_n$ given
by $z_n(t):=\frac{u_n}{\rho_n}$, $\rho_n:=1+\|u_n(0)\|_{H^1}$,  is a $T$-periodic
solution of
\begin{equation}\label{20102015-0735}
\frac{\part z}{\part t} =  \lma_n \Delta z+ \lma_n \rho_n^{-1} f(t,x,\rho_n z), \, x\in\R^N,
\, t>0.
\end{equation}
It is also clear that $v_n:[0,+\infty)\to H^1(\R^N)$ given by $v_n(t): = z_n (t/\lma_n)$ satisfies
\begin{equation}\label{20102015-0736}
\frac{\part v}{\part t} = \Delta v+ \rho_n^{-1} f(t/\lma_n, x,\rho_n v),\, x\in\R^N,\, t>0,
\end{equation}
and that $\rho_n\to +\infty$.
Define $g_n: [0,+\infty) \times \R^N \times \R\to\R$,  $g_n(t,x,v):= \rho_n^{-1}
f(t/\lma_n, x,\rho_n v)$, $n\geq 1$, $t\geq 0$, $x\in\R^N$, $v
\in\R$.
Since the functions $g_n$, $n\geq 1$, satisfy (\ref{23082013-1208}) with a common $m_0$ and (\ref{23082013-1100}) with common $l$ and $\{v_n(0)\}_{n\geq 1}$ is bounded, by use  of Lemma \ref{09012014-1600} and Remark \ref{22092015-1343}, we obtain a constant $R_0>0$ such that $\|v_n(t)\|_{H^1}\leq R_0$ for all $n\geq 1$ and $t\geq 0$.
For a moment fix  an arbitrary $M>0$ and for any $n\geq 1$ take an integer $k_n\geq 1$ such that  $k_n\lma_n T > M$. Observe that Lemma
\ref{26082013-2204} gives, for all $m \geq 1$ and $n\geq 1$,
$$
\|(1-\chi_m) v_n(0)  \|_{L^2}^2 = \|(1-\chi_m) v_n (k_n \lma_n
T)\|_{L^2}^2 \leq R_0^2 e^{-2a k_n \lma_n T} +\alpha_m \leq R_0^2
e^{-2aM}+ \alpha_m
$$
with $\alpha_m \to 0^+$ as $m\to+\infty$. Since $M>0$ is arbitrary we see that $\|(1-\chi_m)  v_n(0)\|_{L^2}\leq \sqrt{\alpha_m}$ for $m, n\geq 1$. Due to the Rellich-Kondrachov for any $m \geq 1$, the set $\{\chi_m v_n(0)\}_{n\geq 1}$ is relatively compact in $L^2(\R^N)$. Therefore $\{ v_n(0) \}_{n\geq 1}$ is relatively compact in $L^2(\R^N)$, since $\alpha_m \to 0^+$ as $m\to +\infty$. As a bounded sequence in $H^1(\R^N)$, $(v_n(0))$ contains a subsequence convergent in $L^2(\R^N)$ to some  $\bar v_0\in H^1(\R^N)$. Therefore we may assume that $v_n(0)\to \bar v_0$ in $L^2(\R^N)$. Moreover, we  may suppose that $\lma_n\to
\lma_0$, as $n\to+\infty$ for some $\lma_0 \in [0,1]$.\\
\indent First consider the case when $\lma_0 \in (0,1]$. Let $f_n: \R^N\times \R\to\R$, $n\geq 1$, be given by
$$
f_n(t,x,z):=\rho_n^{-1}  f(t,x,\rho_n z ), \ \ \mbox{ for all } t\geq 0, \ z\in\R \mbox{ and a.e. } x\in\R.
$$
Note that  (\ref{03092013-0751}) and
(\ref{23082013-1208}) yield
$$
\lim_{n\to+\infty} f_n (t,x,z )= \omega (t,x) z, \mbox{ for all
} t\geq 0, \, z\in\R \mbox{ and  a.e. } x\in\R,
$$
and $\|f_n(t,\cdot,0)\|_{L^2}=\rho_n^{-1} \|f (t,\cdot, 0)\|_{L^2} \leq \rho_n^{-1} \|m_0\|_{L^2}\to
0$, as $n \to +\infty$. It allows us to apply Proposition
\ref{28082013-0033} to (\ref{20102015-0736}). As a result we infer that $z_n (t) \to z_0 (
t)$ in $H^1(\R^N)$ uniformly with respect to $t$ from compact subsets of $(0,+\infty)$, where $z_0:[0,+\infty)\to H^1(\R^N)$ is a $T$-periodic solution of
$$
\frac{\part z}{\part t} = \lma_0 \Delta z + \lma_0 \omega (t,x)z.
$$
Since $\|z(0)\|_{H^1}=\| v(0)\|_{H^1}\neq 0$,  we get a nontrivial $T$-periodic solution of (\ref{11102013-0626}) with $\lma=\lma_0$, a contradiction proving the desired assertion.\\
\indent In the  situation when $\lma_0=0$, we apply Proposition \ref{28082013-0034} to (\ref{20102015-0736}) to see that $v_n(t)\to \widehat v(t)$ uniformly with respect to $t$ from compact subsets of $[0,+\infty)$, where  $\widehat v:[0,+\infty)\to H^1(\R^N)$ is a nontrivial solution of
$$
\frac{\part v}{\part t} = \Delta v + \widehat \omega (x)v, \ x\in \R^N, \ t>0.
$$
Now observe that, for any $t>0$ and $k_n:=[t/\lma_n T]$, $n\geq 1$,  one has
$$
v_n(0)=v_n(k_n\lma_n T) \to \widehat v (t)  \ \mbox{ in }
H^1(\R^N), \mbox{ as } n\to+\infty.
$$
Hence $\widehat v \equiv \bar v_0$ and, as a consequence,
$$
0= \Delta \bar v_0(x) + \widehat \omega (x) \bar v_0 (x), \,
x\in\R^N,
$$
which contradicts the assumption and completes the proof of (i).\\
\indent  To see (ii), suppose that assertion does not hold. Then there exist $\lma_n\in
(0,1)$ and a $T$-periodic solutions $u_n:[0, +\infty)\to H^1(\R^N)$
of
$$
\frac{\part u}{\part t} = \lma_n \Delta u + \lma_n f(t,x,u), \,
x\in\R^N, \, t>0
$$
with $\| u_n (0)\|_{H^1}>0$, $n\geq 1$, and $\| u_n (0)\|_{H^1}\to
0^+$ as $n \to +\infty$. Put $z_n:=\frac{u_n}{\rho_n}$ and let
$v_n(t):=z_n (t/\lma_n)$ with $\rho_n:=\|u_n(0)\|_{H^1}$. Then, for each
$n \geq 1$, $z_n$ is a solution of
$$
\frac{\part z}{\part t} = \lma_n \Delta z+ \lma_n \rho_n^{-1} f(t, x,\rho_n
z),\, x\in\R^N,\, t>0.
$$
and $v_n$ is a solution of
$$
\frac{\part v}{\part t} =  \Delta v+ \rho_n^{-1} f(t/\lma_n, x,\rho_n
v),\, x\in\R^N,\, t>0.
$$
The rest of the proof goes along the lines of the proof for (i). \hfill $\square$

\begin{remark} \label{05092014-1148} $\mbox{ }$\\
{\em  Let us remark that the nonexistence of solutions for
(\ref{11102013-0626}) or (\ref{11102013-0627}) may be also verified
if $\alpha$ or $\omega$ are time dependent. Assume that
\begin{equation}\label{01102014-2200}
\sup_{t \in [0,T]}\|\omega_0(\cdot,t)\|_{L^p} < \left\{
\begin{array}{lcll}
\frac{p^{1/2p}\bar \omega_\infty^{1-1/2p}}{2^{1/2p}},&&
\textnormal{ if } N=1,\ p>2,
\\
\frac{p^{1/p}\bar \omega_\infty^{(1-1/p)}}{4^{1/p}},&&
\textnormal{ if } N=2,\ p>2,
\\
\frac{\bar\omega_\infty^{1-N/2p}}{(N/2p)^{N/2p} C(N)^{N/p}},&&
\textnormal{ if } N\geq 3,  N\leq p<\infty,
\end{array}  \right.
\end{equation}
\noindent where $C(N)>0$ is the constant in the Sobolev inequality
$\|u\|_{L^{\frac{2N}{N-2}}}\leq C(N) \|\nabla u\|_{L^2}$, $u\in
H^1(\R^N)$. Suppose that $u$ is a nonzero $T$-periodic solution of
(\ref{11102013-0626}). Then , for all $t>0$,
\begin{equation}
\frac{\d}{\d t} \frac{1}{2\lambda} \|u(t)\|_{L^2}^2 = -
\int_{\R^N}|\nabla u(t)|^2\d x -  \int_{\R^N} \omega_\infty
(t,x)|u(t)|^2 \d x + \int_{\R^N}\omega_0(t,x)|u(t)|^2 \d
x.\label{04022015-0059}
\end{equation}
Assume first that $N=1$. Then, by use of the H\"{o}lder inequality,
\begin{align}
\int_{0}^{T}\left(\|\nabla u(t)\|_{L^2}^2 +\bar \omega_\infty
\|u(t)\|_{L^2}^2 \right)\d t &\leq \int_0^T
\|\omega_0(t,\cdot)\|_{L^p}\|u(t)\|_{L^2}^{2-2/p}\|u(t)\|_{L^\infty}^{2/p}\d
t, \nonumber \\
&\leq 2^{1/p} \int_0^T \|\omega_0(t,\cdot)\|_{L^p}\|\nabla u(t)\|_{L^2}^{1/p}\|u(t)\|_{L^2}^{2-1/p}\d t,\nonumber%\label{22012015-2253}
\end{align}
where the latter inequality follows by the fact that
$\|u\|_{L^\infty}^2\leq 2 \|\nabla u\|_{L^2}\|u\|_{L^2}$ for $u
\in H^1(\R).$ In the Young inequality $a b
\leq\frac{a^r}{\epsilon^r r} + \frac{b^s \epsilon^s}{s}$ where
$a,b\geq 0$, $\epsilon >0$ and $r\in (1,+\infty)$ such that
$\frac{1}{r}+\frac{1}{s}=1$, put $a:=\|\omega_0 (\cdot,t)\|_{L^p}
\|\nabla u(t)\|_{L^2}^{1/p}$, $b:= \|u(t)\|_{L^2}^{2-1/p}$ and
$r:=2p$ to obtain
\begin{equation}
\|\omega_0(t,\cdot)\|_{L^p}\|\nabla
u(t)\|_{L^2}^{1/p}\|u(t)\|_{L^2}^{2-1/p} \leq
\frac{\|\omega_0(t,\cdot)\|_{L^p}^{2p}\|\nabla
u(t)\|_{L^2}^2}{2p\epsilon^{2p}}+\frac{\epsilon^{2p/(2p-1)}\|u(t)\|_{L^2}^2}{2p/(2p-1)}
\nonumber %\label{22012015-2254}
\end{equation}
for any $\epsilon>0$ and fixed $t \in [0,T]$. If we take
$\epsilon=\epsilon(t)$ so that
$2^{1/p}\frac{\|\omega_0(t,\cdot)\|_{L^p}^{2p}}{2p\epsilon^{2p}}=1$,
i.e. $
\epsilon(t):=\bigg(\frac{2^{(1-p)/p}}{p}\bigg)^{1/2p}\|\omega_0(t,\cdot)\|_{L^p}
$ and apply (\ref{01102014-2200}), then
\begin{align}
\bar\omega_\infty \int_0^T \|u(t)\|_{L^2}^2\d t &\leq
2^{1/p}\frac{2p-1}{2p}\int_0^T
\epsilon(t)^{2p/(2p-1)}\|u(t)\|^2_{L^2}\d t\nonumber \\
& \leq 2^{1/p} \bigg(\frac{2^{(1-p)/p}}{p}\bigg)^{1/(2p-1)} \sup_{t \in [0,T]}\|\omega_0(t,\cdot)\|_{L^p}^{2p/(2p-1)}\int_0^T \|u(t)\|_{L^2}^2\d t \nonumber \\
& = \bigg(\frac{2}{p}\bigg)^{1/(2p-1)} \sup_{t \in [0,T]}\|\omega_0(t,\cdot)\|_{L^p}^{2p/(2p-1)}\int_0^T \|u(t)\|_{L^2}^2\d t \nonumber \\
&<\bar\omega_\infty \int_0^T \|u(t)\|_{L^2}^2\d t,\nonumber
\end{align}
a contradiction proving that (\ref{11102013-0626}) has no nontrivial
$T$-periodic solutions. Assume now that $N=2$. Then by
(\ref{04022015-0059}) and the H\"{o}lder inequality it follows that
\begin{eqnarray*} \nonumber
\int_{0}^{T}\!\!\left(\|\nabla u(t)\|_{L^2}^2 +\bar \omega_\infty
\|u(t)\|_{L^2}^2 \right)\d t \!\leq\!
\int_{0}^{T}\!\! \|\omega_0(t,\cdot)\|_{L^p} \|u(t)\|_{L^{4}}^{4/p}\|u(t)\|_{L^2}^{2-4/p} \d t,\\
\end{eqnarray*}
which in view of the Sobolev inequality $\|u\|_{L^{4}}^2\leq
2\|u\|_{L^2}\|\nabla u\|_{L^2}$, $u \in H^1(\R^2)$ implies
\begin{align}
\int_{0}^{T}\!\!\left(\|\nabla u(t)\|_{L^2}^2 +\bar \omega_\infty
\|u(t)\|_{L^2}^2 \right)\!\d t \!&\leq\!2^{2/p}\!\int_0^T
\|\omega_0(t,\cdot)\|_{L^p} \|\nabla
u(t)\|_{L^2}^{2/p}\|u(t)\|_{L^2}^{2-2/p} \d t.
\label{23012015-0017}
\end{align}
By use of the Young inequality, we obtain for any $\epsilon>0$ and
fixed $t\in [0,T]$,
\begin{align}
\|\omega_0(t,\cdot)\|_{L^p}\|\nabla
u(t)\|_{L^2}^{2/p}\|u(t)\|_{L^2}^{2-2/p}\leq
\frac{\|\omega_0(t,\cdot)\|_{L^p}^p \|\nabla
u(t)\|_{L^2}^2}{p\epsilon^p}+
\frac{\epsilon^{p/(p-1)}\|u(t)\|_{L^2}^2}{p/(p-1)}.\label{23012015-0018}
\end{align}
Choose $\epsilon=\epsilon(t)>0$ such that $\frac{2^{2/p}}{p
\epsilon^p}\|\omega_0(\cdot,t)\|_{L^p}^p=1$, i.e.
$\epsilon(t):=\frac{2^{2/p^2}}{p^{1/p}}\|\omega_0(\cdot,t)\|_{L^p}.
$ Then, by applying (\ref{01102014-2200}), we have
\begin{align}
\bar \omega_\infty \int_0^T \|u(t)\|_{L^2}^2 & \leq
2^{2/p}\frac{(p-1)}{p} \int_0^T
\epsilon(t)^{p/(p-1)}\|u(t)\|_{L^2}^2 \d t \nonumber \\ &
\leq \bigg(\frac{4}{p}\bigg)^{1/(p-1)}\sup_{t \in [0,T]}\|\omega_0(t,\cdot)\|_{L^p}^{p/(p-1)}\int_0^T \|u(t)\|_{L^2}^2 \d t \nonumber \\
& <\bar \omega_\infty \int_0^T \|u(t)\|_{L^2}^2 \d t,\nonumber
\end{align}
which contradicts the existence of nonzero $T$-periodic solution for
(\ref{11102013-0626}). Finally, for $N \geq 3$, by use of the
H\"{o}lder and the Sobolev inequalities we get
\begin{eqnarray*} \nonumber
\int_{0}^{T}\!\!\left(\|\nabla u(t)\|_{L^2}^2 +\bar \omega_\infty
\|u(t)\|_{L^2}^2 \right)\d t \!\leq\!
\int_{0}^{T}\!\! \|\omega_0(t,\cdot)\|_{L^p} \|u(t)\|_{L^{2N/(N-2)}}^{N/p} \|u(t)\|_{L^2}^{2-N/p}\d t\\
\!\leq\! C(N)^{N/p} \int_{0}^{T} \|\omega_0
(\cdot,t)\|_{L^p}\|\nabla u(t)\|_{L^2}^{N/p}
\|u(t)\|_{L^2}^{2-N/p} \d t.
\end{eqnarray*}
In view of the Young inequality, for any $\epsilon>0$ and fixed
$t\in [0,T]$,
$$
\|\omega_0(\cdot,t)\|_{L^p}\|\nabla u(t)\|_{L^2}^{\frac{N}{p}}
\|u(t)\|_{L^2}^{2-\!\frac{N}{p}}\!\!\! \leq\!
\frac{N/2p}{\epsilon^{\frac{2p}{N}}}
\|\omega_0(\cdot,t)\|_{L^p}^{\frac{2p}{N}} \|\nabla
u(t)\|_{L^2}^{2} +\! \left(\!1\!-\!\frac{N}{2p}\! \right)
\epsilon^{\frac{2p}{2p-N}}\!\|u(t)\|_{L^2}^{2} $$ Take
$\epsilon=\epsilon(t)$ so that $\frac{N}{2 p}
\cdot\frac{C(N)^{N/p}}{\epsilon(t)^{\frac{2p}{N}}}
\|\omega_0(\cdot,t)\|_{L^p}^{\frac{2p}{N}}  =1$, i.e.
$\epsilon(t)=(N/2p)^{N/2p} C(N)^{N^2/2p^2}
\|\omega_0(\cdot,t)\|_{L^p}$ and apply (\ref{01102014-2200}), then
\begin{eqnarray*}
& & \bar \omega_\infty \int_{0}^{T} \|u(t)\|_{L^2}^{2}\d t
\leq  C(N)^{N/p} \left(1-\frac{N}{2p}\right)  \displaystyle{\int_{0}^{T}}\epsilon(t)^{\frac{2p}{2p-N}}\|u(t)\|_{L^2}^{2} \d t \\
& & \ \ \ \ \ \ \ \  \leq
(N/2p)^{N/(2p-N)} C(N)^{2N/(2p-N)} \displaystyle{\int_{0}^{T}}\|\omega_0(\cdot,t)\|_{L^p}^{\frac{2p}{2p-N}}\|u(t)\|_{L^2}^{2} \d t\\
& & \ \ \ \ \ \ \ \ \leq (N/2p)^{N/(2p-N)} C(N)^{2N/(2p-N)} \sup_{t\in [0,T]} \|\omega_0(\cdot,t)\|_{L^p}^{\frac{2p}{2p-N}} \int_{0}^{T} \|u(t)\|_{L^2}^{2}\d t\\
& & \ \ \ \ \ \ \ \ < \bar \omega_\infty \int_{0}^{T}
\|u(t)\|_{L^2}^{2}\d t,
\end{eqnarray*}
a contradiction proving that (\ref{11102013-0626}) has no nontrivial
$T$-periodic solutions. }\end{remark}

\section{Proofs of Theorems \ref{02092013-0753} and \ref{06092013-1209}}

We start with a linearization method for computing the fixed point
index of the translation operator in the autonomous case.
\begin{proposition}\label{02092013-0755}
Assume that $f:\R^N\times \R \to\R$ satisfies conditions
{\em(\ref{23082013-1208})}, {\em (\ref{23082013-1100})} and
{\em(\ref{23082013-1209})} {\em(}in their time-independent
versions{\em)} and let ${\bold \Phi}_t$ be the translation along
trajectories for the autonomous equation
$$
\left\{
\begin{array}{l}
\displaystyle{\frac{\part u}{\part t}}(x,t)=
\Delta u(x,t) + f(x,u(x,t)), \ x\in \R^N, \, t>0,\\
u(\cdot, t)\in H^1(\R^N), \, t\geq 0.
\end{array}
\right.
$$
\noindent {\em (i)} If {\em (\ref{03092013-0751})} holds and
$\mathrm{Ker} (\Delta + \omega) =\{ 0\}$, then there exists $R_0>0$
such that  $-\Delta u(x)=f(x,u(x)),\ x\in\R^N$, has no solutions
 $u\in H^1(\R^N)$ with $\|u\|_{H^1} \geq R_0$ and there exists
$\bar t>0$ such that, for all $t\in (0, \bar t]$, ${\bold \Phi}_t
(\bar u)\neq \bar u$ for all $\bar u\in H^1(\R^N) \setminus
B_{H^1}(0,R_0)$ and, for all and all $t\in (0,\bar t]$ and $R\geq
R_0$,
$$
\Ind ({\bold \Phi}_t, B_{H^1} (0,R) ) = (-1)^{m(\infty)}
$$
where $m(\infty)$ is the total multiplicity of the positive eigenvalues  of $\Delta+\omega$.\\
{\em (ii)} If  {\em (\ref{06092013-1208})} holds and $\mathrm{Ker}
(\Delta +\alpha) =\{ 0\}$, then there exists $r_0>0$ such that
$-\Delta u(x)=f(x,u(x)),\ x\in\R^N$, has no solutions
 with $0<\|u\|_{H^1} \leq r_0$ and there exists
$\bar t>0$ such that, for all $t\in (0, \bar t]$, ${\bold \Phi}_t
(\bar u)\neq \bar u$ for all $\bar u\in B_{H^1}(0,r_0)\setminus
\{0\}$ and, for each  $t\in (0,\bar t]$,
$$
\Ind ({\bold \Phi}_t, B_{H^1} (0,r_0) ) =  (-1)^{m(0)}
$$
where $m(0)$ is the total multiplicity of the positive
eigenvalues  of $\Delta+\alpha$.
\end{proposition}

\begin{remark}\label{03092014-1732} {\em
Recall the known arguments on the spectrum of $-\Delta -\omega_0
+\omega_\infty$. To this end, define  ${\bold B}_0: D({\bold
B}_0)\to L^2(\R^N)$ with $D({\bold B}_0):=H^1(\R^N)$ by ${\bold B}_0
u:= \omega_0 u$ and ${\bold B}_\infty: L^2(\R^N)\to L^2(\R^N)$ by
${\bold B}_\infty u:=\omega_\infty u$. By \cite{Pazy}, ${\bold
A}_\infty:={\bold A} - {\bold B}_0 + {\bold B}_\infty$ is  a $C_0$
semigroup generator and its spectrum $\sigma({\bold A}_\infty)$ is
contained in an interval $(-c, +\infty)$ with some $c>0$. It is
clear that $\sigma({\bold A}+{\bold B}_\infty) \subset [\bar
\omega_\infty,+\infty)$. Moreover, it is known, that ${\bold B}_0
({\bold A} + {\bold B}_\infty)^{-1}:L^2(\R^N)\to L^2(\R^N)$ is a
compact linear operator -- for the proof we refer to \cite[Lem.
3.1]{Prizzi-FM}, where the result is obtained under assumption $N
\geq 3$. However a proper restatement, i.e. exploiting Sobolev
embeddings $H^1(\R) \subset L^\infty(\R)$ for $N=1$ and $H^1(\R^2)
\subset L^{4}(\R^2)$ - in case $N=2$ together with the
Rellich-Kondrachov Theorem, leads to the same conclusion. Therefore,
by use of the Weyl theorem
%Since ${\bold B}_0 ({\bold A} + {\bold
%B}_\infty)^{-1}:L^2(\R^N)\to L^2(\R^N)$ is a compact linear operator
%-- see \cite[Lem. 3.1]{Prizzi-FM}, by use of the Weyl theorem
on essential spectra, we obtain $\sigma_{ess} ({\bold A}_\infty)
=\sigma_{ess} ({\bold A}+{\bold B_\infty}) \subset \sigma ({\bold
A}+{\bold B_\infty})\subset [\bar \omega_\infty, +\infty)$ (see e.g.
\cite{Schechter}). Hence, by general characterizations of essential
spectrum, we see that   $\sigma({\bold A}_\infty)\cap(-\infty,0)$
consists of isolated eigenvalues with finite dimensional eigenspaces
(see \cite{Schechter}). }

\end{remark}

\noindent {\bf Proof of Proposition \ref{02092013-0755}:} (i) We
start with an observation that there exists $R_0>0$ such that the
problem
\begin{equation} \label{15102013-0211}
0=\Delta u + (1-\mu)f(x,u) + \mu \omega(x)u, \, x\in\R^N,
\end{equation}
has no weak solutions in $H^1(\R^N)\setminus B_{H^1}(0,R_0)$. To
see this, suppose to the contrary that there exist a sequence
$(\mu_n)$ in $[0,1]$ and solutions $\bar u_n$, $n\geq 1$,  of
(\ref{15102013-0211}) with $\mu=\mu_n$ such that $\|\bar
u_n\|_{H^1}\to +\infty$ as $n\to +\infty$. Put $\rho_n:=1+\|\bar
u_n\|_{H^1}$ and observe that $\bar v_n:= \frac{\bar u_n}{\rho_n}$
are solutions of
$$0=\Delta v + (1-\mu_n) \rho_n^{-1} f(x,\rho_n v)+\mu_n\omega(x)v, \, x\in\R^N.$$
Clearly
$$
\rho_n^{-1} f(x,\rho_{n} v)\to \omega(x)v \mbox{ as } n \to
+\infty \mbox{ for all } t\geq 0 \mbox{ and a.a. } x\in\R^N.
$$
Hence, by use of Remark \ref{15102013-0250} we see that $(\bar u_n)$
contains a sequence convergent to some $\bar u_0\in H^1(\R^N)$ being
a weak nonzero solution of $0=\Delta u + \omega(x) u, \,x\in\R^N,$
a contradiction proving that (\ref{15102013-0211}) has no solutions outside some ball $B_{H^1}(0,R_0)$.\\
\indent Now consider the equation
\begin{equation}\label{02092013-1654}
\frac{\part u}{\part t} = \Delta u + (1-\mu ) f(x,u) +\mu\omega(x)
u, \, x\in\R^N, \, t>0,
\end{equation}
where $\mu\in [0,1]$ is a parameter. Let ${\bold \Psi}_t :
H^1(\R^N)\times [0,1] \to H^1(\R^N)$, $t>0$, be the parameterized
translation along trajectories operator for the above equation. In
view of Theorem \ref{26082013-2118}, there exists $\bar t>0$ such
that
$$
{\bold \Psi}_t (\bar u, \mu)\neq \bar u \ \ \mbox{ for all }\ \ t\in
(0, \bar t], \bar u\in \part B_{H^1}(0,R_0).
$$
By Proposition \ref{01092013-1435} (iii), the homotopy ${\bold
\Psi}_t$ is admissible in the sense of the fixed point theory for
ultimately compact maps (see Section 2). Therefore using the
homotopy invariance one has, for $t\in (0,\bar t]$,
\begin{equation}\label{14112013-1245}
\Ind ({\bold \Phi}_t, B_{H^1}(0,R_0)) = \Ind (e^{-t{\bold
A}_\infty}, B_{H^1}(0,R_0))
\end{equation}
where ${\bold A}_\infty:={\bold A} -{\bold B}_0 + {\bold B}_\infty$.\\
\indent It is left to determine the fixed point index of
$e^{-t{\bold A}_\infty}$.  We note that the set $\sigma({\bold
A}_\infty) \cap (-\infty,0)$ is bounded and closed. Hence, in view
of the spectral theorem (see \cite{Taylor}) there are closed
subspaces $X_-$ and $X^+$ of $L^2 (\R^N)$ such that $X_-\oplus
X^+=L^2 (\R^N)$, $\dim X_-<+\infty$, ${\bold A}_\infty (X_-)\subset
X_-$, ${\bold A}_\infty (D({\bold A_\infty})\cap X^+)\subset X^+$,
$\sigma({\bold A}_\infty|_{X_-})=\sigma({\bold A}_\infty)\cap
(-\infty,0)$, $\sigma({\bold A}_\infty|_{X^+})=\sigma({\bold
A}_\infty)\cap (0,+\infty)$. Define $ {\bold
\Theta}_t:H^1(\R^N)\times [0,1] \to H^1(\R^N)$ by
$$
{\bold \Theta}_t (\bar u, \mu):= (1-\mu) e^{-t {\bold A}_\infty}\bar
u + \mu e^{-t {\bold A}_\infty} {\bold P}_-\bar u,
$$
where ${\bold P}_-:H^1(\R^N)\to H^1(\R^N)$ is the restriction of the
projection onto $X_-\cap H^1(\R^N)$ in $L^2(\R^N)$. Since $\dim
X_-<+\infty$ we infer that ${\bold P}_-$ is continuous. W also claim
that ${\bold \Theta}_t$ is ultimately compact. To see this take a
bounded set $V\subset H^1(\R^N)$ such that $V =
\overline{\conv}^{H^1} {\bold \Theta}_t ( V\times [0,1]).$ This
means that $V \subset \overline{\conv}^{H^1} e^{-t{\bold
A}_\infty}(V\cup {\bold P}_- V)$. Since $V\cup {\bold P}_-V$ is
bounded, Proposition \ref{01092013-1435} (ii) implies that $V$ is
relatively compact in $H^1(\R^N)$, which proves the ultimate
compactness of ${\bold \Theta}_t$. Since $\mathrm{Ker} (I - {\bold
\Theta}_t (\cdot,\mu))=\{ 0\}$ for $\mu\in  [0,1]$, by the homotopy
invariance and the restriction property of the Leray-Schauder fixed
point index, one gets
\begin{eqnarray*}
\Ind (e^{-t{\bold A}_\infty}, B_{H^1}(0,R_0) ) & =&
\Ind_{LS}(e^{-t {\bold A}_\infty}{{\bold P}_-} , B_{H^1}(0,R_0) )\\
& =& \Ind_{LS} ( e^{-t ({\bold A}_\infty|_{X_-})} ,
B_{H^1}(0,R_0)\cap X_-) = (-1)^{m(\infty)}.
\end{eqnarray*}
The latter equality comes from the fact that $\sigma ({\bold
A}_{\infty}|_{X_-})\subset (-\infty,0) $ consists of isolated
eigenvalues of finite dimensional eigenspaces.
This ends the proof of (i) together with (\ref{14112013-1245}).\\
\indent (ii) First we shall prove the existence of $r_0 >0$ such
that the problem
\begin{equation} \label{24012014-0047}
0=\Delta u + (1-\mu)f(x,u) + \mu \alpha(x)u, \, x\in\R^N,
\end{equation}
has no solutions in $B_{H^1}(0, r_0)\setminus \{0\}$. Suppose to the
contrary that there exist a sequence  $(\mu_n)$ in $[0,1]$ and
solutions $\bar u_n:[0,+\infty)\to H^1(\R^N)$, $n\geq 1$,  of
(\ref{24012014-0047}) with $\mu=\mu_n$ such that $\|\bar
u_n\|_{H^1}\to 0^+$ as $n\to +\infty$ and $\|\bar u_n\|_{H^1}\neq
0$, $n\geq 1$. Put $\rho_n:=\|\bar u_n\|_{H^1}$. Then $\bar
v_n:=\frac{\bar u_n}{\rho_n}$ are solutions of
$$0=\Delta v + (1-\mu_n) \rho_n^{-1} f(x,\rho_n v)+\mu_n\alpha(x)v, \, x\in\R^N.$$
Observe that
$$
\rho_n^{-1} f(x,\rho_{n} v)\to \alpha(x)v \mbox{ as }  n \to \infty
\mbox{ for a.a. } x\in\R^N.
$$
Using again Remark \ref{15102013-0250} one can see that $(\bar u_n)$
(up to a subsequence) converges to some nonzero solution of
$0=\Delta u + \alpha (x) u, \,x\in\R^N,$ a contradiction. Summing
up, there is $r_0 > 0$ such that (\ref{24012014-0047}) has no
solutions $u \in H^1(\R^N)$ with $0<\|u\|_{H^1} \leq r_0$. The rest
of the proof runs as before: by ${\bold \Psi}_t : H^1(\R^N)\times
[0,1] \to H^1(\R^N)$, $t>0$ we denote the translation along
trajectories operator for the equation
\begin{equation}\label{24012014-0112}
\frac{\part u}{\part t} = \Delta u + (1-\mu ) f(x,u) +\mu\alpha(x)
u, \, x\in\R^N, \, t>0, \, \mu\in [0,1],
\end{equation}
and, by applying Theorem \ref{26082013-2118} we obtain the existence
of $\bar t>0$ such that
$$
{\bold \Psi}_t (\bar u, \mu)\neq \bar u \ \ \mbox{ for all }\ \ t\in
(0, \bar t], \bar u\in \part B_{H^1}(0,r_0).$$ Next Proposition
\ref{01092013-1435} (iii) ensures the admissibility of ${\bold
\Psi}_t$ and by homotopy invariance, for $t\in (0,\bar t]$, we have
\begin{equation}\label{24012014-0119} \Ind ({\bold \Phi}_t,
B_{H^1}(0,r_0)) = \Ind (e^{-t{\bold A}_0}, B_{H^1}(0,r_0))
\end{equation}
where ${\bold A}_0:={\bold A}-{\bold C}_0 + {\bold C}_\infty$ and
operators ${\bold C}_i: D({\bold C}_i) \to L^2(\R^N)$ with $D({\bold
C}_i) = H^1(\R^N)$ are given by ${\bold C}_i u:= \alpha_i u$, $i \in
\{0,\infty\}$. Now one can easily determine fixed point index of
$e^{-t{\bold A}_0}$ by arguing as in part (i) (with ${\bold
A}_\infty$ replaced by ${\bold A}_0$ and $B_{H^1}(0, R_0)$ replaced
by $B_{H^1}(0, r_0)$) and, as a consequence, obtain that
\begin{eqnarray*}
\Ind (e^{-t{\bold A}_0}, B_{H^1}(0,r_0) ) & =&
\Ind_{LS}(e^{-t {\bold A}_0}{{\bold P}_-} , B_{H^1}(0,r_0) )\\
& =& \Ind_{LS} ( e^{-t ({\bold A}_0|_{X_-})} , B_{H^1}(0,r_0)\cap
X_-) = (-1)^{m (0)}.
\end{eqnarray*}
This completes the proof.
\hfill $\square$\\

\vspace{5mm}

\noindent Now we are ready to conclude and provide proofs of our main results.\\

\noindent {\bf Proof of Theorem \ref{02092013-0753}}: Let ${\bold
\Phi}_t$, $t>0$, be the translation operator for
(\ref{23082013-0153}). It is clear that
$$
\lim_{|u|\to +\infty}\frac{\widehat f(x,u)}{u} = \widehat \omega(x),
\,  \mbox{ for any } x\in\R^N.
$$
Hence, by applying Proposition \ref{02092013-0755} (i) we obtain
$R_0 > 0$ such that
$$
\Delta u(x) + \widehat f(x,u(x))=0, \ x\in\R^N,
$$
has no solutions in the set $H^1(\R^N)\setminus B_{H^1}(0,R_0)$ and
there exists $t_0>0$ such that, for $t\in (0,t_0]$,
\begin{equation}\label{141002013-1344}
\Ind (\widehat {\bold \Phi}_t, B_{H^1}(0,R_0)) = (-1)^{m(\infty)}.
\end{equation}
Due to Proposition \ref{04102013-2013} and the assumption,
increasing $R_0$ if necessary, we can assume that
(\ref{11102013-0626}) has no $T$-periodic solutions starting  from
$H^1(\R^N) \setminus B_{H^1}(0,R_0)$. Taking $U:=B_{H^1}(0,R_0)$ and
applying Corollary \ref{14102013-1246} we get
$$
\Ind ({\bold \Phi}_T, B_{H^1}(0,R_0))= \lim_{t\to 0^+} \Ind
(\widehat{\bold \Phi}_t, B_{H^1}(0,R_0)),
$$
which along with (\ref{141002013-1344}) yields $\Ind ({\bold
\Phi}_T, B_{H^1}(0,R_0))=(-1)^{m(\infty)}$.
This and the existence property of the fixed point index imply that there exists $\bar u\in B_{H^1}(0,R_0)$ such that ${\bold \Phi}_T (\bar u)=\bar u$, i.e. there exists a $T$-periodic solution of (\ref{23082013-0153}).  \hfill $\square$\\

\noindent {\bf Proof of Theorem \ref{06092013-1209}}: First use
Proposition \ref{02092013-0755} to get $R_0,r_0 >0$ such that
\begin{equation}\label{02012014-2107}
\lim_{t\to 0^+} \Ind (\widehat {\bold \Phi}_t, B_{H^1}(0,R)) =
(-1)^{m_-(\infty)} \ \mbox{ if } \ R\geq R_0
\end{equation}
and
\begin{equation}\label{02012014-2108}
\lim_{t\to 0^+} \Ind (\widehat {\bold \Phi}_t, B_{H^1}(0,r)) =
(-1)^{m(0)} \ \mbox{ if } \ 0 < r\leq r_0.
\end{equation}
Now, due to Proposition \ref{04102013-2013} there exist $R\geq R_0$
and $r\in (0,r_0]$ such that, for any $\lma\in (0,1]$,
(\ref{09102013-1413}) has no solutions with $u(0) \in B_{H^1}(0,r)
\cup \left( H^1(\R^N) \setminus B_{H^1}(0,R)\right)$. Next we put
$U:= B_{H^1}(0,R) \setminus \overline{B_{H^1}(0,r)}$ and apply
Corollary \ref{14102013-1246} to get
$$
\Ind({\bold\Phi}_T, U) = \lim_{t\to 0^+} \Ind (\widehat {\bold
\Phi}_t, U).
$$
This together with  (\ref{02012014-2107}) and (\ref{02012014-2108}),
by use of the additivity property of the fixed point index, yields
\begin{eqnarray*}
\Ind({\bold\Phi}_T, U) & = & \lim_{t\to 0^+} \Ind (\widehat {\bold
\Phi}_t, B_{H^1}(0,R)) -
\lim_{t\to 0^+} \Ind (\widehat {\bold \Phi}_t, B_{H^1}(0,r))\\
& = & (-1)^{m(\infty)}-(-1)^{m(0)} \neq 0,
\end{eqnarray*}
which gives the existence of the fixed point of ${\bold \Phi}_T$ in $U$.  \hfill $\square$\\

% \begin{example}
% {\em  (a) Consider the following problem
% $$
% \left\{ \begin{array}{l}
% \displaystyle{\frac{\part u}{\part t}} = \Delta u + % \left(\mu+\frac{1}{|x|^s}+\sin t\right) u +  g(\lma u), \ % \ x\in\R^3, \, t>0,\\
% u(\cdot,t)\in H^1(\R^3), \ t\geq 0,\\ u(x,t)=u(x, t+2\pi), x\in \R^3,\, t\geq 0,
% \end{array} \right. $$
% where $\mu:[0,\infty)\to \R$, $s>1$,

% (b) Let $f:[0,+\infty)\times \R^N\times \R\to\R$ be given by
% $$
% f(t,x,u):=a_\infty (x,t)u+a_0 (x,t)u + g(b_\infty (x)u)) + h(b_0(x)u) + k(t,x), \, t\geq
% $$

% $$  \frac{\part u}{\part t} = \Delta u + , \, x\in\R^N, \, t>0, $$
% where $f, g:\R\to\R$ and $f'(0)\neq 0$, $g'(0)\neq 0$.
% If $f(t,x,u):=(a(x,t)u+b(x,t))u+g(a_0(x)u)+h(b_0(x)u)$.
% }\end{example}


\begin{thebibliography}{99}


\bibitem{Akhmerov-et-al}
R.R. Akhmerov, M.I. Kamenskii , A.S. Potapov, A.E. Rodkina, B.N.
Sadovskii, {\em Measures of Noncompactness and Condensing
Operators}, Birkh\"{a}user 1992.


\bibitem{Amman-Zehnder} H. Amman, E. Zehnder , {\em Nontrivial solutions for a class of nonresonance problems and applications to nonlinear differential equations}, Annali della Scuola Superiore di Pisa, t. 7 (4), 1980, 539--603.

\bibitem{Antoci-Prizzi} F. Antoci, M. Prizzi, {\em  Attractors and global averaging of non-autonomous reaction-diffusion equations in $\R^N$}, Topol. Methods Nonlinear Anal. 20 (2002), no. 2, 229--259.

\bibitem{Becker} R.R. Becker, {\em Periodic solutions of semilinear equations of evolution of compact type},
J. Math. Anal. Appl. 82 (1981) 33--48.

\bibitem{BreNir} H. Brezis, L. Nirenberg, {\em Characterizations of the ranges of some nonlinear operators
and applications to boundary value problems}, Ann. Scuola Norm. Sup.
Pisa Cl. Sci. (4) 5 (1978), no. 2, 225--326.

\bibitem{Cholewa} J. Cholewa, T. D{\l}otko, \emph{Global Attractors in Abstract Parabolic
Problems}, Cambridge University Press, 2000.


\bibitem{Cwiszewski-CEJM} A. \'Cwiszewski, {\em Positive periodic solutions of parabolic evolution problems: a translation
along trajectories approach}, Centr. Eur. J.  Math., vol. 9, no. 2
(2011), 244--268.

\bibitem{Cwiszewski-hab} A. \'Cwiszewski, {\em Periodic and stationary solutions of nonlinear evolution
equations: translation along trajectories method}, habilitation
report, Toru\'n, 2011, http://arxiv.org/pdf/1309.6295.pdf.

\bibitem{Cw-Luk-2015} A. \'{C}wiszewski, R. {\L}ukasiak, {\em A Landesman�Lazer type result for periodic parabolic problems on $R^N$ at resonance}, Nonlinear An. TMA 125 (2015), 608--625.

% \bibitem{Daners} D. Daners, P. Koch Medina, \emph{Abstract evolution equations, periodic problems and applications}, Longman Scientific $\&$ Technical, New York 1992.

\bibitem{Deimling} K. Deimling, {\em Multivalued Differential Equations}, Walter de Gruyter, Berlin, New York, 1992.

\bibitem{Dunford-Schwartz} N. Dunford, J. T. Schwartz,
{\em Linear Operators}, Parts I and II, Wiley-Interscience, New York
1966.

% \bibitem{Evans} L. C. Evans, \emph{Partial Differential Equations},
% Graduate Studies in Mathematics 19, American Mathematical Society, Providence, RI, 1998.

% \bibitem{Trudinger} D.Gilbarg, N.S. Trudinger, \emph{Elliptic Partial Differential Equations of Second Order}, Springer, 1997.

\bibitem{Granas} A. Granas, J. Dugundji, {\em Fixed Point Theory}, Springer-Verlag, New York 2003.

\bibitem{Henry} D. Henry, \emph{Geometric Theory of Semilinear
Parabolic Equations}, Springer Verlag, 1981.

\bibitem{Hess} P. Hess, {\em Periodic-parabolic boundary value problems and positivity}, Pitman Research
Notes in Mathematics Series, 247, Longman Scientific \& Technical,
John Wiley \& Sons, 1991.

\bibitem{Hirano1} N. Hirano, {\em Existence of multiple periodic solutions for a semilinear evolution equations},
Proc. Amer. Math. Soc. 106 (1989), 107-114.

\bibitem{Hirano2} N. Hirano, {\em Existence of periodic solutions for nonlinear evolution equations in
Hilbert spaces}, Proc. Amer. Math. Soc. 120, no. 1 (1994), 107--114.

\bibitem{Hu-Papa} S. Hu, N. S. Papageorgiou N., {\em On the existence of periodic solutions for a class of nonlinear
evolution inclusions}, Boll. Unione Mat. Ital., 7B, 1993, 591--605.

\bibitem{NkaWill} M. N. Nkashama, M. Willem, {\em Time-periodic solutions of boundary value problems
for nonlinear heat, telegraph and beam equations}, in Differential
Equations: Qualitative Theory, Szeged (Hungary); Coll. Math. Soc.
J�nos Bolyai, 1984, 47, 809--845.

\bibitem{Pazy} A. Pazy, \emph{Semigroups of linear operators and
applications to partial differential equations}, Springer Verlag
1983.

\bibitem{Prizzi-FM} M. Prizzi, {\em On admissibility of parabolic equations in $\R^N$}, Fund. Math. 176 (2003), 261--275.

\bibitem{Pruss} J. Pr\"{u}ss, {\em Periodic solutions of semilinear evolution equations}, Nonlinear Anal. 3,
(1979), 221--235.

\bibitem{Schechter} M. Schechter, {\em Spectra of partial differential operators}, North-Holland 1986.

\bibitem{Sell-You} G. Sell, Y. You, {\em Dynamics of Evolutionary Equations}, Springer Verlag, 2002.

\bibitem{Shioji} N. Shioji, {\em Existence of periodic solutions for nonlinear evolution equations with nonmonotonic perturbations}, Proc. Amer. Math. Soc., 125 (1997), 2921--2929.

\bibitem{Tanabe} H. Tanabe, \emph{Equations of evolution}, Monographs and Studies in Mathematis no. 6, Pitman 1979.

\bibitem{Taylor} A. E. Taylor, {\em  Functional Analysis}, Wiley, New York 1961.

\bibitem{Vejvoda} O. Vejvoda, L. Herrmann, V. Lovicar, M. Sova, I. Stra\v{s}kraba, M. \v{S}t\v{e}dr\'{y}, {\em Partial differential equations: time-periodic solutions}, Martinus Nijhoff Publishers, The Hague, 1981.

\bibitem{Wang} B. Wang, {\em Attractors for reaction-diffusion equations in unbounded domains}, Phys. D 128 (1999), 41--52.

\end{thebibliography}
\end{document}